\documentclass[11pt,twoside]{article}
\usepackage{epsfig,amsfonts,color}
\usepackage{amsmath,bm}
\usepackage{mathrsfs}
\allowdisplaybreaks
\usepackage{multirow}
\usepackage{cite}
\usepackage{amssymb, palatino, geometry,url}

\usepackage[colorlinks=true,linkcolor=blue,citecolor=blue,urlcolor=blue]{hyperref}
\usepackage{lineno}
\usepackage{enumerate}
\usepackage{float}
\usepackage{graphicx}
\usepackage{array}
\usepackage{placeins}
\usepackage{caption}
\usepackage{subcaption}
\usepackage{pdflscape}
\usepackage{algorithm}
\usepackage{algpseudocode}
\usepackage{fancyhdr}
\newcommand{\mods}[1]{\left|  #1 \right|}

\newcommand{\be}{\begin{equation}}
\newcommand{\ee}{\end{equation}}
\newcommand{\bea}{\begin{eqnarray}}
\newcommand{\eea}{\end{eqnarray}}

\newcommand{\bvec}{\left(\begin{array}{c}}
\newcommand{\evec}{\end{array}\right)}
\newcommand{\bsub}{\begin{subequations}}
\newcommand{\esub}{\end{subequations}}

\newcolumntype{L}[1]{>{\raggedright\let\newline\\\arraybackslash\hspace{0pt}}m{#1}}
\newcolumntype{C}[1]{>{\centering\let\newline\\\arraybackslash\hspace{0pt}}m{#1}}
\newcolumntype{R}[1]{>{\raggedleft\let\newline\\\arraybackslash\hspace{0pt}}m{#1}}

\begin{document}

\title{Dual Dynamic Programming for \\ Multi-Scale Mixed-Integer MPC}

\author{Ranjeet Kumar${}^{\P}$, Michael J. Wenzel${}^\ddag$, Mohammad N. ElBsat${}^\ddag$, Michael J. Risbeck${}^\ddag$,\\ Kirk H. Drees${}^\ddag$, Victor M. Zavala${}^{\P}$\\
{\small ${}^{\P}$Department of Chemical and Biological Engineering}\\
{\small \;University of Wisconsin-Madison, 1415 Engineering Dr, Madison, WI 53706, USA}\\
{\small ${}^\ddag$Johnson Controls International}\\
{\small \;507 E. Michigan St., Milwaukee, WI 53202, USA}}
 \date{}
\maketitle

\begin{abstract}                
We propose a dual dynamic integer programming (DDIP) framework for solving multi-scale mixed-integer model predictive control (MPC) problems. Such problems arise in applications that involve long horizons and/or fine temporal discretizations as well as mixed-integer states and controls (e.g., scheduling logic and discrete actuators). The approach uses a nested cutting-plane scheme that performs forward and backward sweeps along the time horizon to adaptively approximate cost-to-go functions. The DDIP scheme proposed can handle general MPC formulations with mixed-integer controls and states and can perform forward-backward sweeps over block time partitions. We demonstrate the performance of the proposed scheme by solving mixed-integer MPC problems that arise in the scheduling of central heating, ventilation, and air-conditioning (HVAC) plants. We show that the proposed scheme is scalable and dramatically outperforms state-of-the-art mixed-integer solvers. 
\end{abstract}


\section{Introduction} \label{sec:intro}

Model predictive control (MPC) is a flexible optimization-based control technology that can handle complex decision-making settings through the use of different types of objectives, variables, and constraints. Of particular recent interest are MPC formulations that can capture  {\em mixed-integer logic} that involves discrete actuators and states (e.g., on/off switches and scheduling logic)  and continuous actuators and states (e.g., inventory and production levels) \cite{Garcia1989,Mayne2000,rawlingsbook,kirches2011fast}. This interest has also been sparked by recent results that show that stability and robustness properties of MPC are applicable to settings that involve discrete decisions  \cite{rawlings2017model,bemporad1999control,di2014stabilizing,aguilera2013stability,rawlings2017model}. Unfortunately, mixed-integer MPC formulations are difficult to handle computationally,  as they involve combinatorial decision spaces that grow with the time horizon and resolution.  This hinders practical industrial applications, as settings of interest often involve decisions that span multiple timescales (e.g., scheduling,  economic optimization, and control). 

An important application of mixed-integer MPC is that of scheduling and control of energy systems;  for instance, the operation of central heating, ventilation, and air-conditioning (HVAC) plants for campuses or urban districts involve continuous and discrete decisions that need to be made hourly and over horizons that span weeks  \cite{rawlings2018,risbeck2015cost}.  Specifically, long horizons and fine time resolutions are required to capture peak demand charges and short-term fluctuations of energy loads and prices \cite{Ranjeet2018,kumar2019,kumar2020hvac}. Large-scale mixed-integer MPC problems are also commonly encountered in integrated production  scheduling and control for manufacturing systems  \cite{daoutidis2018integrating,baldea2015integrated,beal2017combined,beal2018integrated}. Chemical production scheduling is another common application that involves mixed-integer optimization formulations and it is well-known than such formulations can be posed in state-space form   \cite{subramanian2012state,gupta2017general,gupta2016deterministic,pattison2017moving,gupta2016rescheduling}.

Advances in algorithms and solvers for mixed-integer linear programming (MILP) have enabled the solution of sophisticated mixed-integer MPC formulations but a wide range of applications are still inaccessible to off-the-shelf tools.  Complexity is often handled via hierarchical decomposition; for instance, high-level (infrequent) scheduling decisions are often decoupled from low-level (frequent) economic optimization and control decisions.  Complexity can also be handle by decomposing the problem over the time domain. This temporal decomposition approach has been used in production scheduling problems \cite{bassett1996decomposition,harjunkoski2001decomposition,jackson2003temporal}, unit commitment in power systems \cite{kim2018temporal}, and hydrothermal scheduling \cite{sifuentes2007hydrothermal}. Existing decomposition paradigms are primarily based on Benders or Lagrangian schemes. An issue with Benders decomposition is that it cannot directly incorporate discrete variables in subproblems; as a result, this approach has been primarily used for hierarchical decomposition (by capturing all discrete logic in the master problem).  Lagrangian decomposition is a more flexible approach that can handle mixed continuous and discrete variables but this approach often exhibits slow convergence. 

In pioneering work, Pereira and co-workers \cite{Pereira1989,Pereira1991} proposed a scalable temporal decomposition scheme to tackle multi-stage stochastic programming problems with purely continuous variables.  This approach is now called stochastic dual dynamic programming (SDDP) and uses nested Benders decomposition scheme that solves a sequence of small subproblems (defined over a single stage or timestep) in forward and backward sweeps along the time horizon. These sweeps update state information (in forward mode) and collect dual information (in backward mode) and this information is used to construct cutting planes that approximate the so-called cost-to-go function (also known as the recourse function in stochastic programming). A key benefit of SDDP is that it scales to problems with many stages (long horizons and/or fine time resolutions);  for instance, we have recently used this approach to tackle multi-stage stochastic programming problems arising in battery management with over 168 stages (such problems are intractable with any other known technique) \cite{kumar2018stochastic}.  Recent work has shown that SDDP can also be used to tackle multi-stage stochastic programs with binary state variables \cite{zou2016,zou2019}. This scheme has been called stochastic dual dynamic integer programming (SDDIP). The authors in \cite{lara2018,lara2019} have used SDDIP to tackle electrical power infrastructure planning problems that involve both continuous and integer variables.  Unfortunately, existing SDDIP schemes cannot be directly applied to mixed-integer MPC formulations of interest because the system is not represented in standard state-space form. 

In this work, we present a dual dynamic integer programming (DDIP) algorithm for solving mixed-integer MPC formulations. The scheme is derived for a general state-space representations that involve both continuous and discrete control and state variables. The ability to handle general state-space representations enables the use of the algorithm in diverse applications such as MPC with discrete actuators and hybrid logic, real-time scheduling, integrated scheduling, and control. The proposed approach can also perform forward-backward sweeps over block time partitions and this enables more effective use of efficient off-the-shelf solvers (compared to a standard scheme that operates over individual timesteps). The block representation can also be used to more naturally capture problems with multiple timescales (a block represents a subscale) and reveals connections between forward-backward sweeps and receding-horizon schemes. Specifically, a receding horizon scheme is a forward sweep that ignores the cost-to-go function (value of future information). The proposed scheme also reveals important connections between dual (adjoint) variables, cost-to-go functions, and value of future information. We demonstrate the performance of the DDIP scheme through computational experiments for a mixed-integer MPC problem for a central HVAC plant for a typical university campus. We show that the DDIP scheme is scalable and can dramatically reduce computational times of state-of-the-art MILP solvers.

\section{Mixed-Integer MPC Formulation} \label{sec:formulation}
We consider a mixed-integer MPC problem of the following form: 
\begin{subequations}\label{eq:mpcproblem}
\begin{align}
\min_{\mathbf{x}_{N},\mathbf{u}_{N}} \; &\sum_{t=0}^{N-1} f_{t}(x_t, u_t) \label{eq:obj}\\
\textrm{s.t.} \; & x_{t+1} = A_tx_{t} + B_tu_{t},\; t = 0,  \dots, N-1  \label{eq:dynamics}\\
& x_0 = \bar{x} \label{eq:initial} \\
& u_{t} \in \mathcal{U}_t, u_t \in \mathbb{Z}^{n_{u_1}} \times \mathbb{R}^{n_{u_2}},  \;  t = 0,  \dots , N-1 \label{eq:controlbounds}\\
& x_{t} \in \mathcal{X}_t, x_t \in \mathbb{Z}^{n_{x_1}} \times \mathbb{R}^{n_{x_2}},  \;  t = 0,  \dots , N.  \label{eq:statebounds}
\end{align}
\end{subequations}
Here, the controls and states at time $t$ are denoted as $u_t$ and $x_t$, respectively. The function $f_t(x_t,u_t)$ is the stage cost at time $t$. The control trajectory over the entire time domain is denoted as $\mathbf{u}_N$ and the state trajectory is denoted as $\mathbf{x}_N$. The initial state is given by $x_0=\bar{x}$. We note that the states $x_t$ and controls $u_t$ are mixed-integer; specifically, these have $n_{u_1}$ and $n_{x_1}$ integer components and $n_{u_2}$ and $n_{x_2}$ continuous components, respectively. In other words, we have that $x_{t} \in \mathbb{Z}^{n_{x_1}} \times \mathbb{R}^{n_{x_2}}$ and $u_{t} \in \mathbb{Z}^{n_{u_1}} \times \mathbb{R}^{n_{u_2}}$ (with $n_x=n_{x_1}+n_{x_2}$, $n_u=n_{u_1}+n_{u_2}$). The states and controls are also bounded by time-varying polyhedral sets $\mathcal{X}_t$ and $\mathcal{U}_t$, respectively. The system dynamics are expressed in state-space form by using the time-varying matrices $A_t \in \mathbb{R}^{n_x \times n_x}$, $B_t \in \mathbb{R}^{n_x \times n_u}$. In this work, we assume that the stage costs $f_t(x_t,u_t)$ are  linear functions of the controls and the states. The above formulation can easily be modified to capture disturbances (these are not expressed explicitly in order to simplify the presentation). 

\subsection{Dynamic Programming Representation} \label{subsec:dynamicprog}

We now reformulate the mixed-integer MPC problem \eqref{eq:mpcproblem} using an equivalent dynamic programming representation \cite{bertsekas1995dynamic}. In this representation, the objective function is expressed as the sum of the current stage cost and a cost-to-go function that embeds the cost of future stages. The MPC problem in dynamic programming format is written as:
\begin{subequations}\label{eq:Q0mpcproblem}
\begin{align}
Q_0(x_0) = \min_{x_{1}, u_{0}} \; & f_{0}(x_0, u_0) + Q_1(x_1) \label{eq:Q0obj}\\
 \textrm{s.t.} \; & x_{1} = A_{0}x_{0} + B_{0}u_{0} \label{eq:Q0dynamics0}\\
& u_{0} \in \mathcal{U}_{0}, u_0 \in \mathbb{Z}^{n_{u_1}} \times \mathbb{R}^{n_{u_2}} \label{eq:Q0controlbounds}\\
& x_{1} \in \mathcal{X}_{1}, x_1 \in \mathbb{Z}^{n_{x_1}} \times \mathbb{R}^{n_{x_2}}.  \label{eq:Q0statebounds}
\end{align}
\end{subequations}
where the state $x_0$ is given as the initial state $\bar{x}$, and the cost-to-go $Q_1(x_1)$ is given by:
\begin{subequations}\label{eq:Q1mpcproblem}
\begin{align}
Q_1(x_1) = \min_{x_{2}, u_{1}} \; & f_{1}(x_1, u_1) + Q_2(x_2) \label{eq:Q1obj}\\
\textrm{s.t.} \; & x_{2} = A_{1}x_{1} + B_{1}u_{1} \label{eq:Q1dynamics0}\\
& u_{1} \in \mathcal{U}_{1}, u_1 \in \mathbb{Z}^{n_{u_1}} \times \mathbb{R}^{n_{u_2}} \label{eq:Q1controlbounds}\\
& x_{2} \in \mathcal{X}_{2}, x_2 \in \mathbb{Z}^{n_{x_1}} \times \mathbb{R}^{n_{x_2}}.  \label{eq:Q1statebounds}
\end{align}
\end{subequations}
Following the recursion, the cost-to-go for stage $t$ is  given by:
\begin{subequations}\label{eq:Qtmpcproblem}
\begin{align}
Q_t(x_t) = \min_{x_{t+1}, u_{t}} \; & f_{t}(x_t, u_t) + Q_{t+1}(x_{t+1}) \label{eq:Qtobj}\\
\textrm{ s.t.} \; & x_{t+1} = A_{t}x_{t} + B_{t}u_{t} \label{eq:Qtdynamics0}\\
& u_{t} \in \mathcal{U}_{t}, u_t \in \mathbb{Z}^{n_{u_1}} \times \mathbb{R}^{n_{u_2}} \label{eq:Qtcontrolbounds}\\
& x_{t+1} \in \mathcal{X}_{t+1}, x_{t+1} \in \mathbb{Z}^{n_{x_1}} \times \mathbb{R}^{n_{x_2}}.  \label{eq:Qtstatebounds}
\end{align}
\end{subequations}
The cost-to-go at the final stage $N-1$ is:
\begin{subequations}\label{eq:Q2mpcproblem}
\begin{align}
Q_{N-1}(x_{N-1}) = \min_{x_{N}, u_{N-1}} \; & f_{N-1}(x_{N-1}, u_{N-1}) \label{eq:QN1obj}\\
\textrm{ s.t.} \; & x_{N} = A_{N-1}x_{N-1} + B_{N-1}u_{N-1} \label{eq:QN1dynamics0} \\
& u_{N-1} \in \mathcal{U}_{N-1}, u_{N-1} \in \mathbb{Z}^{n_{u_1}} \times \mathbb{R}^{n_{u_2}} \label{eq:QN1statebounds}\\
& x_{N} \in \mathcal{X}_{N}, x_{N} \in \mathbb{Z}^{n_{x_1}} \times \mathbb{R}^{n_{x_2}}. \label{eq:QN1controlbounds}
\end{align}
\end{subequations}

The dynamic programming representation of the MPC problem reveals that future information is embedded in the cost-to-go function but such function does not have a closed form. Our goal will be to derive schemes that approximate the cost-to-go functions. To enable this, we introduce an auxiliary variable $z_t \in \mathbb{R}^{n_x}$ (as a copy of $x_t$) in Problem \eqref{eq:Qtmpcproblem}. We note that the auxiliary variable $z_t$ is a continuous variable but satisfies all the constraints of the mixed-integer state $x_t$ at the solution. 
\begin{subequations}\label{eq:Qtauxproblem}
\begin{align}
Q_t(x_t) = \min_{x_{t+1}, u_{t},z_{t}} \; & f_{t}(x_t, u_t) + Q_{t+1}(x_{t+1}) \label{eq:Qtauxobj}\\
\textrm{ s.t.} \; & x_{t+1} = A_{t}z_{t} + B_{t}u_{t} \label{eq:Qtauxdynamics0}\\
& z_{t} = x_{t}  \quad (\mu_{t}) \label{eq:Qtauxvariable}\\ 
& u_{t} \in \mathcal{U}_{t}, u_t \in \mathbb{Z}^{n_{u_1}} \times \mathbb{R}^{n_{u_2}} \label{eq:Qtauxcontrolbounds}\\
& x_{t+1} \in \mathcal{X}_{t+1}, x_{t+1} \in \mathbb{Z}^{n_{x_1}} \times \mathbb{R}^{n_{x_2}}  \label{eq:Qtauxstatebounds}\\
& z_{t} \in \mathbb{R}^{n_x}.  \label{eq:Qtauxvariablebounds}
\end{align}
\end{subequations}
The introduction of the  auxiliary variable $z_t$ induces a lifting procedure that more clearly reveals the coupling structure of the MPC problem. Specifically, note that $z_t$ is the only variable that couples the cost-to-functions between stages (via the linking constraint \eqref{eq:Qtauxvariable}). As such, if the linking constraint is removed, the initial state of the stage suproblem becomes a free variable (thus decoupling the problem). The dual variable (adjoint) for the linking constraint \eqref{eq:Qtauxvariable} in the corresponding linear programming (LP) relaxation of Problem \eqref{eq:Qtauxproblem} is denoted as $\mu_t$. This dual variable represents the value of information of the state $x_t$ and will be key in designing the proposed dual dynamic integer programming (DDIP) scheme. 

\subsection{Dual Dynamic Integer Programming} \label{subsec:dualdynamic}
We use the dynamic programming representation \eqref{eq:Qtauxproblem} to design the DDIP scheme. The proposed scheme adapts the SDDIP framework presented in \cite{zou2016,zou2019} to handle deterministic MPC problems in general state-space form.  In the proposed scheme, the cost-to-go functions are approximated by using a collection of cutting planes (hyperplanes) accumulated in each iteration. However, cutting planes can only be constructed if the cost-to-go functions are continuous and convex. As such, we obtain the cutting planes from linear programming (LP) relaxations of the stage subproblems. 

Each iteration $k = 0,1,\dots$ of the DDIP algorithm comprises a forward sweep and a backward sweep along the time horizon. The forward sweep in iteration $k$ marches sequentially in time to compute a control trajectory $(u_{0}^{k},u_{1}^{k},\dots,u_{N-1}^{k})$ and state trajectory $(x_{1}^{k},x_{2}^{k},\dots,x_{N}^{k})$ by solving the stage subproblems by using given approximate cost-to-go functions $\phi_{t}^{k}(\cdot)$ (instead of the actual functions $Q_{t}(\cdot)$). The approximate cost-to-go $\phi_{t}^{k}(\cdot)$ is given by Problem \eqref{eq:Qtapp}. In the backward sweep of iteration $k$, a set of Benders cuts  (one for each stage $t=N-1,N-2,\dots,0$) is generated by solving the LP relaxation of Problem \eqref{eq:Qtapp} using given values for the states obtained in the forward sweep. The solution to the LP relaxation of the Problem \eqref{eq:Qtapp} gives updated values for the cost-to-go $\hat{\phi}_k^t$ and for the dual variable $\mu_{t}^{k}$ for the linking constraint, which are used to construct Benders cuts that approximate the cost-to-go functions. 

\begin{subequations}\label{eq:Qtapp}
\begin{align}
\phi_{t}^{k}(x_{t}^{k}) = \min_{x_{t+1}, u_{t}, z_{t},\theta_{t+1}} \; & f_{t}(x_{t}^{k}, u_{t}) + \theta_{t+1} \label{eq:Qtapp_obj}\\
\textrm{s.t.} \; & x_{t+1} = A_{t}z_{t} + B_{t}u_{t} \label{eq:Qtapp_dynamics0} \\
& z_{t} = x_{t}^{k}  \quad (\mu_{t}^{k}) \label{eq:Qtappvariable}\\ 
& \theta_{t+1} \geq \phi_{t+1}^{l} + (\mu_{t+1}^{l})^T (x_{t+1}-x^l_{t+1}), \; l = 0, 1, \dots, k-1. \label{eq:Qtapp_cut}\\  
& u_{t} \in \mathcal{U}_{t}, u_t \in \mathbb{Z}^{n_{u_1}} \times \mathbb{R}^{n_{u_2}} \label{eq:Qtappcontrolbounds}\\
& x_{t+1} \in \mathcal{X}_{t+1}, x_{t+1} \in \mathbb{Z}^{n_{x_1}} \times \mathbb{R}^{n_{x_2}}  \label{eq:Qtappstatebounds}\\
& z_{t} \in \mathbb{R}^{n_x}.  \label{eq:Qtappvariablebounds}
\end{align}
\end{subequations}

The constraints \eqref{eq:Qtapp_cut} are the collection of Benders cuts for stage $t$ that have been accumulated over all previous iterations. In each iteration $k$ and for each stage $t$, the variable $\theta_{t+1}$ is used to approximate the cost-to-go for the next stage. We can justify the choice of the cut coefficients in \eqref{eq:Qtapp_cut} by using LP duality properties. Specifically, for the LP relaxation of Problem \eqref{eq:Qtapp}, if the dual solution corresponding to the linking constraint \eqref{eq:Qtappvariable} is $\mu^k_{t}$, the corresponding cost-to-go can also be represented in dual form as $\hat{\phi}^k_t(x^k_t) = {(\mu^{k}_{t})}^T x^k_{t} + v^k_t$, where $v^k_t$ denotes the cost component in dual form corresponding to the constraints defining the polyhedral sets $\mathcal{U}_{t}$ and $\mathcal{X}_{t}$. We can see that if $\mu^{k}_{t}$ and $v^k_{t}$ are substituted using information from the previous iteration by $\mu^{k-1}_{t}$ and $\hat{\phi}^{k-1}_{t} - {(\mu^{k-1}_{t})}^T x^{k-1}_{t}$, respectively, we obtain the Benders cut for stage $t-1$ for the current iteration $k$: $\theta_t \geq {(\mu^{k-1}_{t})}^T (x^k_{t} - x^{k-1}_{t}) + \hat{\phi}^{k-1}_{t}$. In each iteration $k$, the objective $\hat{\phi}_{t}^{k}(x_{t}^{k})$ and the dual solution $\mu_{t}^{k}$ obtained at stage $t$ in the backward sweep are used to define the Benders cut for the problem at stage $t-1$ for the next iteration $k+1$, which provides a lower bound for the cost-to-go $\phi_{t}^{k}(x_{t}^{k})$ that appears in the objective at stage $t-1$.
\\

In a forward sweep, the stage subproblem minimizes the current stage cost ($f_{t}(x_{t},u_{t})$) and the cost-to-go (approximated by the cutting planes). The state variable $x_{t}$ at time $t$ encodes information of the history of the control actions implemented in the previous stages, while the cost variable $\theta_{t+1}$ encodes information of the future stages. We note that, in iteration $k$, solving Problem \eqref{eq:Qtapp} in the forward sweep gives a set of trajectories for the controls and states $\{u_t^k,x_{t+1}^k\}, t=0,\dots,N-1$, starting with the given initial state $x_0=\bar{x}$. These control and state policies typically satisfy the constraints of Problem \eqref{eq:mpcproblem}; in other words, they constitute a feasible (but suboptimal) solution.  Consequently, the cost obtained from the forward sweep $\sum_{t=0}^{N-1}f_t(x_t^k,u_t^k)$ is  guaranteed to be an upper bound for the MPC objective. We can interpret the first forward sweep of the DDIP scheme as a receding-horizon MPC scheme that computes an approximate policy of the long-horizon problem by ignoring future information. In subsequent iterations, the forward sweep is interpreted as a receding horizon scheme that is guided by the approximate cost-to-go functions. In some cases, the forward sweep might lead to infeasibility (e.g., neglecting future information might put the system at a state that cannot satisfy future constraints). Consequently, the use of approximate cost-to-go functions can be critical in ensuring feasibility. In the proposed framework, infeasibilities are handled by adding slack variables and penalty terms in the cost function. It is also possible to handle infeasibility in DDIP by adding feasibility cuts but this approach is more complicated and difficult to implement. 
\\

In a backward sweep, cutting planes are collected and added to the approximate problems $\phi_t^k(\cdot)$ of each stage which incrementally brings the objective function closer to the optimal. The duals  $\mu_t$ of the LP relaxation problem are time-dependent adjoint variables that propagate information of the future backward in time. We note that the objective value $\hat{\phi}^k_0(x_0^k)$ provides a lower bound for the MPC objective $Q_0(x_0)$ because $\hat{\phi}^k_0(x_0^k)$ is obtained from the solution to the LP relaxation of Problem \eqref{eq:Qtapp} and also contains lower approximations of cost-to-go functions at all stages. The algorithm is stopped when the upper and lower bounds are within a pre-defined tolerance $\epsilon$. 
\\

\newpage 

The DDIP scheme can be summarized as:
\emph{
\begin{enumerate}
\item START with $x_{0}^{0} = \bar{x}$, $\hat{\phi}_{t+1}^{0} = 0$, $\mu_{t+1}^{0} = 0$ for $t=0,1,\dots,N-1$, and set iteration counter $k \leftarrow 0$.\\ \\
\textbf{Forward Pass:}
\item Compute $\phi_{0}^{k}(x_{0}^{k})$, and get $x_{1}^{k}$, $u_{0}^{k}$. 
\item For $t=1,2,\dots,N-1$, solve $\phi_{t}^{k}(x_{t}^{k})$ and get the policy $x_{t+1}^{k}$, $u_{t}^{k}$.
\item Set the current upper bound, \\$ub = \sum\limits_{t=0}^{N-1}{f_t(x_t^k,u_t^k)}$. \\ \\
\textbf{Backward Pass:} 
\item For $t = N-1,N-2,\dots,1$, solve LP relaxation of the updated Problem \eqref{eq:Qtapp}, and obtain $\hat{\phi}_{t}^{k}(x_{t}^{k})$ and associated dual $\mu_{t}^{k}$. 
\item For  $t = N-1,N-2,\dots,1$, compute the corresponding Benders cuts: $\theta_{t} \geq \hat{\phi}_{t}^{k} + {(\mu_t^k)}^T (x_t-x^k_t)$; add the computed cuts to the problems $\phi_{t-1}^{k+1}(x_{t-1}^{k+1})$ for future iterations.
\item Set the current lower bound, $lb = \hat{\phi}_{0}^{k}(x_{0}^{k})$. 
\item UPDATE $k \leftarrow k+1$. If ${ub - lb} \geq \epsilon$, RETURN to Step 2, else STOP.
\end{enumerate}}

A key advantage of the DDIP scheme is that the time horizon is decomposed into stages and it only needs to solve single-timestep optimization problems. In other words, the DDIP algorithm moves sequentially over time and never needs to form and solve the entire mixed-integer MPC problem explicitly. This enables reductions in computational time and overcomes computer memory bottlenecks. The cutting planes that approximate the cost-to-go functions encode (summarize) all the accumulated information. We also note that in each iteration, the backward sweep of the DDIP scheme solves single-stage problems with fixed initial states provided from the forward sweep, and therefore each stage subproblem is independent of other stage subproblems. Consequently, all stage subproblems in the backward sweep can be solved {\em in parallel} and this enables further reductions of computational time. 


\subsection{Convergence Analysis} \label{subsec:convergence}
We note that, at iteration $k$, computing $\phi_{t}^{k}(x_t)$ in the forward sweep gives a set of controls, states, and auxiliary variable $(u_t^k,x_{t+1}^k,z_t^k)$ which satisfy the constraints of the problem $Q_t(x_t)$. Consequently, $(u_t^k,x_{t+1}^k,z_t^k)$ is a feasible solution to $Q_t(x_t)$ (but not necessarily optimal) and the computed value o f$\phi_{t}^{k}(x_t)$ is an upper bound of $Q_{t}(x_t)$. We also note that the backward sweep delivers a lower bound for the objective function if the cuts generated in the backward sweep are valid. If the cut coefficients in iteration $k$, given by $(\hat{\phi}_{t}^{k}, \mu_t^k)$, are obtained from the backward sweep, and the Benders cut $\hat{\phi}_{t}^{k} + {(\mu_t^k)}^T (x_t-x^k_t)$ is generated with the state $x^k_t$ computed in the forward pass in iteration $k$ and fixed for the backward pass, then the cut is said to be valid for stage $t$ and iteration $k$ if it satisfies \cite{zou2019}:
\begin{align}
\hat{\phi}_{t}^{k} + {(\mu_t^k)}^T (x_t-x^k_t) \leq Q_{t-1}(x^k_{t-1}) \label{eq:validcut}
\end{align}
In other words, the cutting plane is valid if it provides an underestimator of the cost-to-go function. Because  the Benders cut $\hat{\phi}_{t}^{k} + {(\mu_t^k)}^T (x_t-x^k_t)$ is generated by solving the LP relaxation of Problem \eqref{eq:Qtapp} with the optimal value $\hat{\phi}_{t}^{k}$, it underestimates the optimal value $\phi_{t}^{k}$ of Problem \eqref{eq:Qtapp} (i.e., $\hat{\phi}_{t}^{k} \leq \phi_{t}^{k}$). It follows that, because the Benders cut is valid for the LP relaxation of Problem \eqref{eq:Qtapp}, it is also valid for the MILP problem \eqref{eq:Qtapp} and therefore, also for the Problem \eqref{eq:Qtauxproblem} \cite{birge1985}.  Unfortunately, convergence of the DDIP scheme cannot be guaranteed because there may be a duality gap associated with the solution (since the problem is mixed-integer and the cost-to-go functions are obtained from LP relaxations).   We will see, however, that this scheme delivers high-quality solutions in practice.  The DDIP algorithm is guaranteed to converge in problems that only involve continuous variables. 


\subsection{Multi-Scale Mixed-integer MPC Scheme} \label{subsec:multiscale}
The established connections between mixed-integer MPC and DDIP reveal possibilities to design more sophisticated decomposition schemes. The DDIP scheme introduced in Section \ref{subsec:dualdynamic} assumes a subproblem defined over a single timestep per stage, and implicitly assumes that information is collected after each timestep. When multiple timescales are present in a problem, however, we can design a DDIP scheme that decomposes the time horizon into stages, in such a way that each stage is a block partition that contains multiple timesteps. This means that the DDIP scheme updates state and dual information at a slower timescale (stages) than the system dynamics (that evolves over inner timesteps). This multiscale design of the DDIP scheme provides computational flexibility to leverage the use of off-the-shelf solvers. Here, we present modifications needed to handle block time partitions. 

Consider the mixed-integer MPC problem in which the DDIP information (cutting planes) are collected at stages $t$ and $t+1$, but the system dynamics evolve over the inner time sequence:
\begin{equation}
t=t_{0}, t_{1}, \dots, t_{n-1}, t_{n}, t_{n+1}=t+1. 
\end{equation}
At stage $t$, given the state $x_t$ at the beginning of the stage (corresponding to inner time $t_0=t$), the cost-to-go $Q_t(x_t)$ is solved to obtain the stage trajectories $\boldsymbol{x_{t+1}} := (x_{t_1},x_{t_2},\dots,x_{t_{n}},x_{t+1})$ and $\boldsymbol{u_{t}} := (u_{t},u_{t_1},u_{t_2},\dots,u_{t_{n}})$. The cost-to-go $Q_t(x_t)$ for stage $t$ is given by the subproblem:
\begin{subequations}\label{eq:Qtapp_multi}
\begin{align}
Q_{t}(x_{t}) = \min_{\boldsymbol{x_{t+1}}, \boldsymbol{u_{t}}} \; & f_t(\boldsymbol{x_{t}},\boldsymbol{u_{t}}) + Q_{t+1}(x_{t+1}) \label{eq:Qtapp_obj_multi}\\
\textrm{s.t.} \; & x_{t_1} = A_{t}z_{t} + B_{t}u_{t} \label{eq:Qtapp_dynamicst1_multi} \\
& x_{{t}_{p+1}} = A_{{t}_p}x_{{t}_p} + B_{{t}_p}u_{{t}_p}, \; p=1,\dots,n. \label{eq:Qtapp_dynamicstp_multi} \\
& z_{t} = x_{t} \quad (\mu_{t}) \label{eq:Qtauxvariable_multi}\\ 
& u_{{t}_p} \in \mathcal{U}_{{t}_p}, u_{{t}_p}  \in \mathbb{Z}^{n_{u_1}} \times \mathbb{R}^{n_{u_2}}, \; p=0,\dots,n. \label{eq:Qtappcontrolbounds_multi}\\
& x_{{t}_p} \in \mathcal{X}_{{t}_p}, x_{{t}_p} \in \mathbb{Z}^{n_{x_1}} \times \mathbb{R}^{n_{x_2}}, \; p=1,\dots,n. \label{eq:Qtappstatebounds_multi}\\
& z_{t} \in \mathbb{R}^{n_x}.  \label{eq:Qtappvariablebounds_multi}
\end{align}
\end{subequations}
From duality, we have that the LP relaxation of this problem satisfies $\hat{Q}_{t}(x_{t}) = (\mu_{t})^T x_{t} + v_t$ similar to the discussion in Section \ref{subsec:dualdynamic}. We thus have that a consistent DDIP scheme can be designed by using a Benders cut obtained for each stage $t$ in the backward sweep in each iteration $k$. We use the computed LP relaxation objective $\hat{\phi}^k_{t}$, state $x^k_{t_1}$, and dual $\mu_t^k$ in iteration $k$ to define the Benders cut for stage $t-1$ for future iterations: $\theta_{t} \geq \hat{\phi}^k_{t} + (\mu_t^k)^T(x_{t_1}-x^k_{t_1})$ (similar to that in \eqref{eq:Qtapp_cut}). We show the implementation of this {\em multi-scale DDIP scheme} in Section \ref{sec:experiments}.

We can solve the subproblem \eqref{eq:Qtapp_multi} with off-the-shelf solvers which can provide good primal solutions very fast if the number of steps included in the subproblems is small. However, we also observe that the above scheme can be used to construct hierarchical MPC schemes that tackle multiple timescales.  For instance, the subproblem \eqref{eq:Qtapp_multi} can be solved using DDIP and one can partition a given stage into multiple substages and solve each substage with DDIP. The implementation of such schemes is an interesting topic of future work. 


\section{Computational Experiments} \label{sec:experiments}
In this section, we present computational experiments that demonstrate the scalability of the proposed DDIP schemes. Here, we use challenging mixed-integer MPC problems that arise in the operation of central HVAC plants. 

\subsection{Decision-Making Setting} \label{subsec:setting}

The decision-making setting of the central HVAC plant studied here is an extension of the continuous formulation described in \cite{kumar2020hvac}. The HVAC plant that we consider consists of a chiller subplant comprising 4 chillers that produce chilled water and a heat recovery (HR) chiller subplant comprising 3 HR chiller units that produce both chilled water and hot water, 3 hot water generators to produce hot water, 9 cooling towers to reduce the temperature of the water purchased from the market, a dump heat exchanger (dump HX) for rejecting heat from the hot water, and storage tanks (one for chilled water and one for hot water). 

The goal of the MPC controller is to determine hourly operating strategies for each equipment unit to satisfy the demands of chilled and hot water from a collection of buildings of a university campus and to minimize the total cost of the external utilities that need to be purchased from the market (electricity, water, and natural gas). Electricity is charged based on time-varying prices, while water and natural gas usage are charged at constant prices. The various cost components for the central plant are:
\begin{itemize}
\item {\em Electricity transactions (hourly)}: Electricity is required for the equipment operation in the central plant. The transactions are charged at the time-varying market price, $\pi^{e}_{t}$. 
\item {\em Water transactions (hourly)}: Water is required to make up for evaporative losses of water in the cooling towers. Water is purchased from the utility at a fixed price $\pi^{w}_{t}$= \$0.009/gal. 
\item {\em Natural gas transactions (hourly)}: Natural gas is needed to run hot water generators to satisfy the campus heating load. Natural gas is purchased from the utility at a fixed price of $\pi^{ng}_{t}$= \$0.018/kWh. 
\end{itemize}
\begin{figure}[!ht]
\centering
\includegraphics[width=\textwidth]{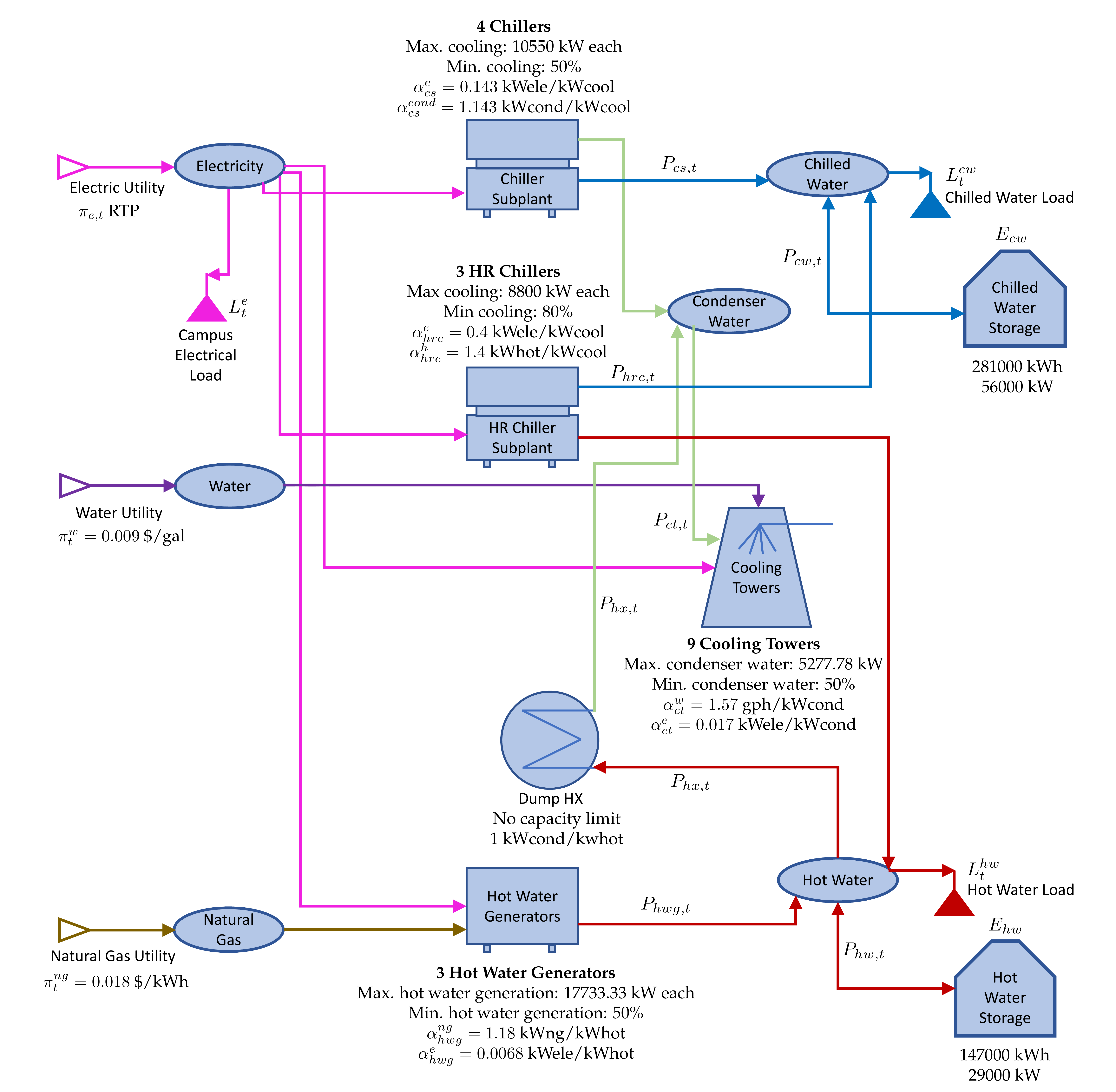}
\caption{Schematic representation of central HVAC plant under study \cite{kumar2020hvac}.}
\label{fig:schematic}
\end{figure}

Figure \ref{fig:schematic} shows the energy flows between all the units of the central HVAC plant and interactions with campus loads and utilities. Electricity is consumed by the chiller subplant, HR chiller subplant, hot water generators, and cooling towers, while utility water is consumed by only the cooling towers to make up for the evaporative losses, and natural gas is only consumed by the hot water generators. The electricity, water, and natural gas consumption of these units are linked to their operating loads. For instance, all chiller units and HR chiller units consume $\alpha^{e}_{cs}$ and $\alpha^{e}_{hrc}$ kW of electricity per kW of chilled water produced, respectively; each hot water generator consumes $\alpha^{e}_{hwg}$ kW of electricity and $\alpha^{ng}_{hwg}$ kW of natural gas per kW hot water produced, respectively; and each cooling tower consumes $\alpha^{e}_{ct}$ kW of electricity and $\alpha^{w}_{ct}$ utility water per kW of condenser water input, respectively. Also, there are prescribed minimum operating loads for the chillers, HR chillers, hot water generators, and cooling towers. For instance, each chiller can either be OFF or operate at 50\% or more of its maximum cooling capacity when ON. Similarly, each HR chiller unit, when ON, can only operate at minimum 80\% of its maximum cooling capacity, each hot water generator and each cooling tower can operate at 50\% or more of the respective maximum capacities, while the dump HX can operate at any load when ON. These restrictions require enforcing logical constraints and making binary decisions for each of these units whether to keep a unit ON or OFF, thus resulting in a mixed-integer formulation for the central plant operations.

The central HVAC plant meets the chilled water load ($L^{cw}_{t}$) of the university campus by producing chilled water from the chiller subplant ($P_{cs,t}$), the HR chiller subplant ($P_{hrc,t}$), and the discharge from chilled water storage ($P_{cw,t}$), where, within the chiller subplant, $n^{th}$ unit produces $p_{n,cs,t}$ chilled water, and within the HR chiller subplant, $n^{th}$ unit produces $p_{n,hrc,t}$ chilled water. The hot water production from the HR chiller subplant ($\alpha^{h}_{hrc}P_{hrc,t}$), the hot water generators ($P_{hwg,t}$), and the discharge from the hot water storage ($P_{hw,t}$) meets the hot water load ($L^{hw}_{t}$) of the campus, where $n^{th}$ HR chiller unit produces $\alpha^{h}_{hrc}p_{n,hrc,t}$ hot water, and $n^{th}$ hot water generator unit produces $p_{n,hwg,t}$ hot water. The excess hot water ($P_{hx,t}$) in the system is recycled by dump heat exchanger (HX) by cooling it and producing condenser water which is cooled further by the cooling towers together with the condenser water produced by the chiller and HR chiller subplants (total $P_{ct,t}$ condenser water is cooled by the towers, where $p_{n,ct,t}$ is cooled by $n^{th}$ tower). The binary control decisions for the system include the ON/OFF decisions for each of the chiller units, HR chiller units, hot water generators, and cooling towers. The continuous control decisions are the operating loads of all units, which include the chilled water production by each chiller and HR chiller unit (when ON), hot water production by each hot water generator (when ON), the cooling load of each cooling tower (when ON), the heat exchange load of the dump HX, and discharge rates from the two storage tanks.

The HVAC plant operations are driven by time-varying disturbances, which are given by the campus loads for electricity ($L^{e}_{t}$), chilled water ($L^{cw}_{t}$), and hot water ($L^{hw}_{t}$), and by the electricity prices ($\pi^{e}_{t}$). In this work, we assume perfect knowledge of these disturbances. 

\subsection{Mixed-Integer MPC Formulation} \label{subsec:mpcproblem}
The mixed-integer MPC formulation uses the disturbance forecasts for $L^{e}_{t}$, $L^{cw}_{t}$, $L^{hw}_{t}$, and $\pi^{e}_{t}$ to find the control and state policy that minimizes the total cost over an $N$-hour prediction horizon, with the time sets $\mathcal{T} := \{0,\dots,N\}$ and $\bar{\mathcal{T}} := \mathcal{T} \setminus \{N\}$, by solving the following mixed-integer linear program: 
\begin{subequations} \label{eq:detmpc}
\begin{align}
\min \; & \sum\limits_{t \in \bar{\mathcal{T}}} \left(\sum\limits_{j=\{e,w,ng\}}\pi^j_t r^{j}_{t} - \pi^e_t L^e_t + \sum\limits_{j \in \{cw,hw\}} \rho_{j}U_{j,t}\right). \label{eq:objective}\\
\textrm{s.t. }\; & Y_{n,j,t} \in \{0,1\}, \;  j \in \{cs,hrc,hwg,ct,hx\}, n \in \{1,\dots,N_j\}, t \in \bar{\mathcal{T}} \label{eq:binary}\\
& \beta_j\overline{p}_jY_{n,j,t} \leq p_{n,j,t} \leq \overline{p}_jY_{n,j,t}, \;  j \in \{cs,hrc,hwg,ct,hx\}, n \in \{1,\dots,N_j\}, t \in \bar{\mathcal{T}} \label{eq:powers} \\
& P_{j,t} = \sum\limits_{n=1}^{N_j}{p_{n,j,t}}, \; j \in \{cs,hrc,hwg,ct,hx\}, t \in \bar{\mathcal{T}} \label{eq:totpowers} \\
& r^{e}_{t} = \sum_{j \in \{cs,hrc,hwg,ct\}} \alpha^{e}_{j} P_{j,t} + L^{e}_{t}, \; t \in \bar{\mathcal{T}} \label{eq:elec_load}\\
& r^{j}_{t} = \alpha^{j}_{i_j} P_{i_j,t}, \; j \in \{w,ng\}, \; i_{w}=ct, i_{ng}=hwg, \; t \in \bar{\mathcal{T}}  \label{eq:water_ng_load} \\
& P_{ct,t} = \alpha^{cond}_{cs} P_{cs,t}+P_{hx,t}, \; t \in \bar{\mathcal{T}} \label{eq:dumphx} \\
& P_{cs,t}+P_{hrc,t}+P_{cw,t}+U_{cw,t}= L^{cw}_{t}, \; t \in \bar{\mathcal{T}} \label{eq:cw_load} \\
& \alpha^{h}_{hrc} P_{hrc,t}+P_{hwg,t}-P_{hx,t}+P_{hw,t}+U_{hw,t}  = L^{hw}_{t}, \; t \in \bar{\mathcal{T}} \label{eq:hw_load} \\
& E_{j,t+1} = E_{j,t} - P_{j,t}, \; j \in \{cw,hw\}, t \in \bar{\mathcal{T}} \label{eq:edynamics} \\
& 0 \leq E_{j,t} \leq \overline{E}_{j,t}, \; j \in \{cw,hw\}, t \in \mathcal{T} \label{eq:ebound}\\
& \underline{P}_j \leq P_{j,t} \leq \overline{P}_j, \; j \in \{cw,hw\}, t \in \bar{\mathcal{T}} \label{eq:psbound}\\
& U_{j,t} \geq 0, \;   j \in \{cw, hw\},  t \in \bar{\mathcal{T}} \label{eq:Ubound} 
\end{align}
\end{subequations}
Here, the binary variables $Y_{n,j,t}$ in \eqref{eq:binary} are used to introduce the logical ON/OFF decisions for each of the chiller, HR chiller, hot water generator, cooling tower, and dump HX units, and constraints \eqref{eq:powers} provide the operating bounds of each of these units when ON. Constraints \eqref{eq:totpowers} define the total operating loads of the chiller subplant, HR chiller subplant, hot water generators, cooling towers, and dump HX, respectively. The constraints \eqref{eq:elec_load}-\eqref{eq:water_ng_load} compute the demands of electricity, water,  and natural gas ($r^{e}_t,r^{w}_t,r^{ng}_t$) that need to be purchased from the utility companies. Constraints \eqref{eq:dumphx} impose the energy balance for the condenser water.  Constraints \eqref{eq:cw_load} and \eqref{eq:hw_load} ensure that the chilled and hot water loads are met. Slack variables $U_{j,k}$, $ j \in \{cw, hw\}$ for unmet loads for chilled and hot water ensure that the problem remains feasible under unmet (under-production) chilled water or hot water demands. The slack variables for unmet loads are penalized in the objective function with a large penalty coefficient $\rho_j$, $j \in \{cw,hw\}$. We will see that demand satisfaction violations can be encountered when the MPC controller does not have a sufficiently long horizon and thus the use of future information (in the form of cost-to-go functions) can help mitigate these issues. 

The dynamics of the state-of-charge (SOC) for chilled and hot water storage tanks are given by constraints \eqref{eq:edynamics}. Constraints \eqref{eq:ebound}-\eqref{eq:Ubound} provide bounds on the states, discharge rates from storage tanks, and slack variables for unmet loads. The lower bounds for the discharge rates of chilled water and hot water storage units correspond to the maximum charging rates, which are negative of the maximum discharging rates (i.e, $\underline{P}_j=-\overline{P}_j, j \in \{cw,hw\}$). The objective function includes cost components for each of the electricity, water, natural gas purchase from the utility, and penalty for the slack variables. We subtract the cost of the baseline electricity load of the campus, $\pi^e_tL^e_t$ from the objective because it is a constant value.

\begin{figure}[!htp]
\centering
\includegraphics[width=0.7\textwidth]{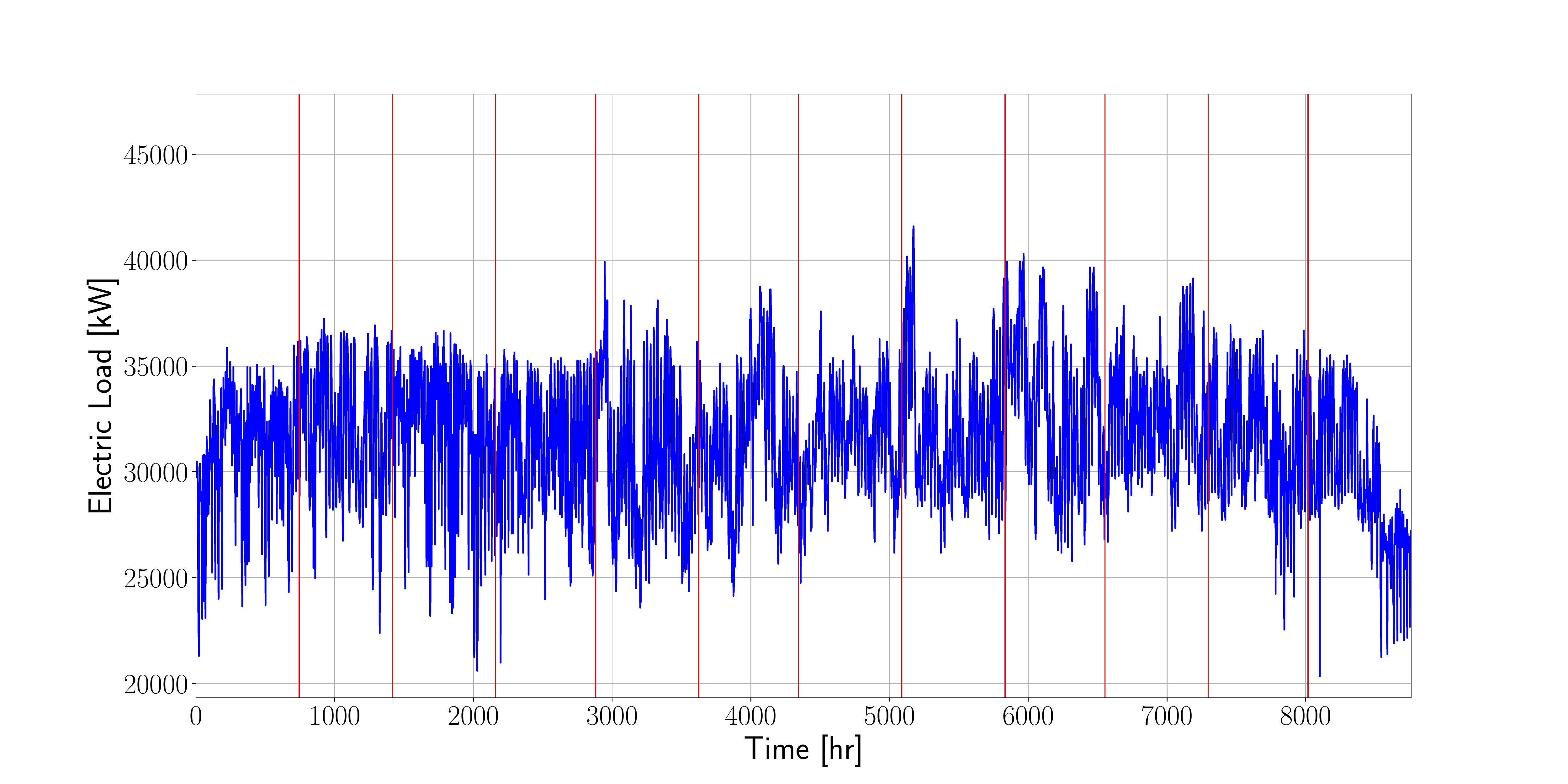}
\caption{Historical electrical load of the campus. Red vertical lines denote end of each month.}
\label{fig:eload}
\end{figure}

\begin{figure}[!htp]
\centering
\includegraphics[width=0.7\textwidth]{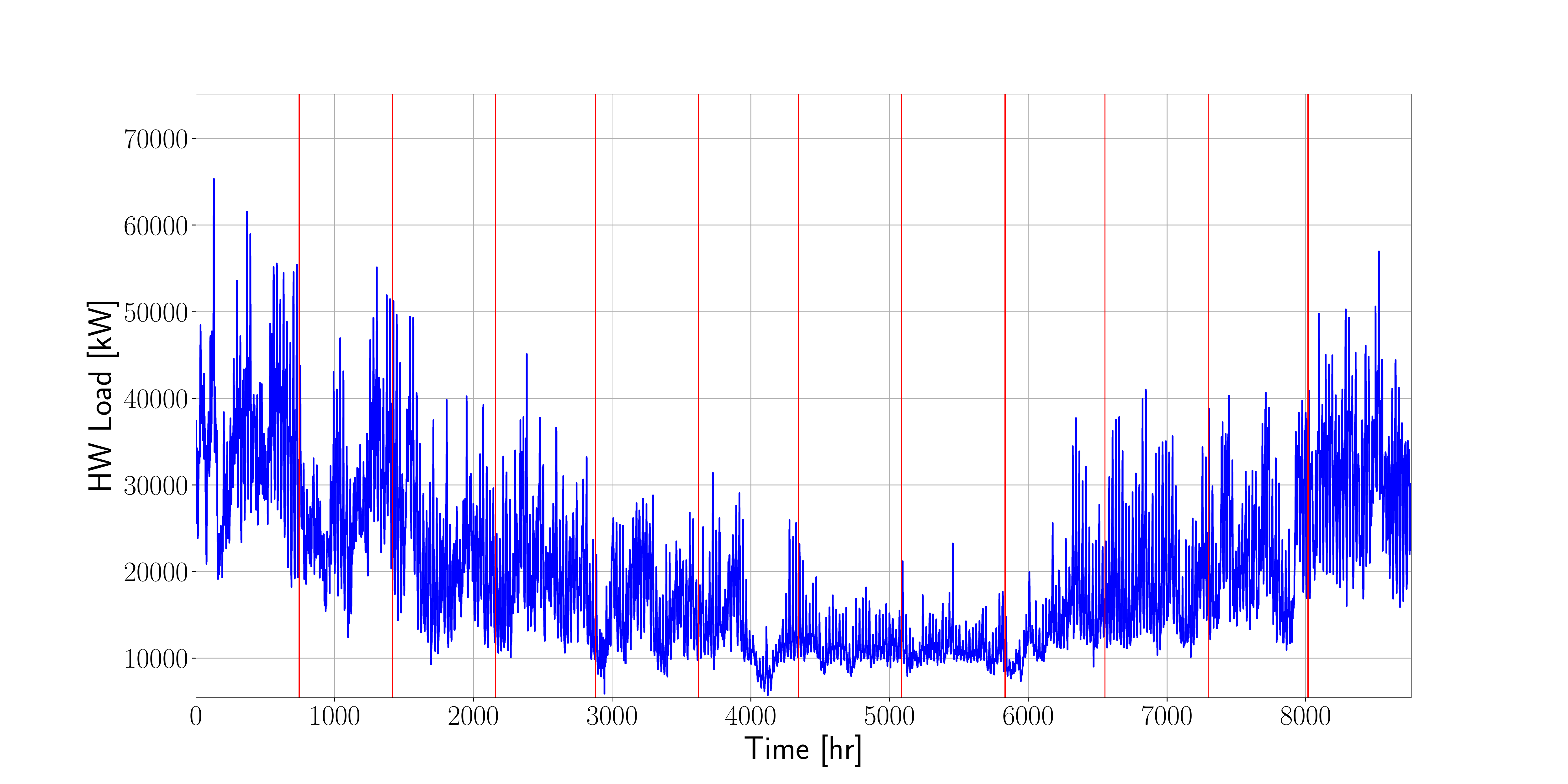}
\caption{Historical hot water load of the campus. Red vertical lines denote end of each month.}
\label{fig:hload}
\end{figure}

\begin{figure}[!htp]
\centering
\includegraphics[width=0.7\textwidth]{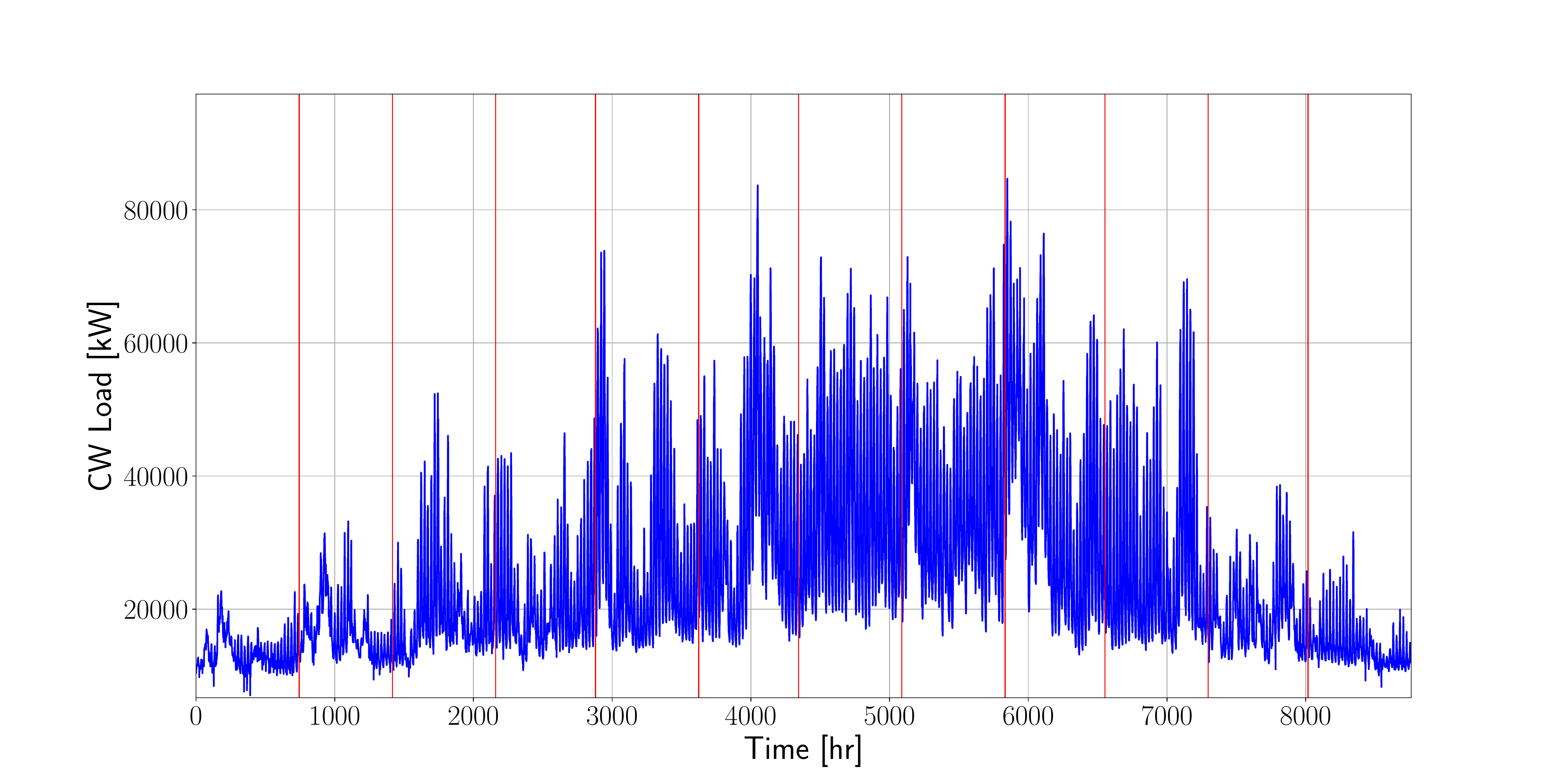}
\caption{Historical chilled water load of the campus. Red vertical lines denote end of each month.}
\label{fig:cload}
\end{figure}

\begin{figure}[!htp]
\centering
\includegraphics[width=0.7\textwidth]{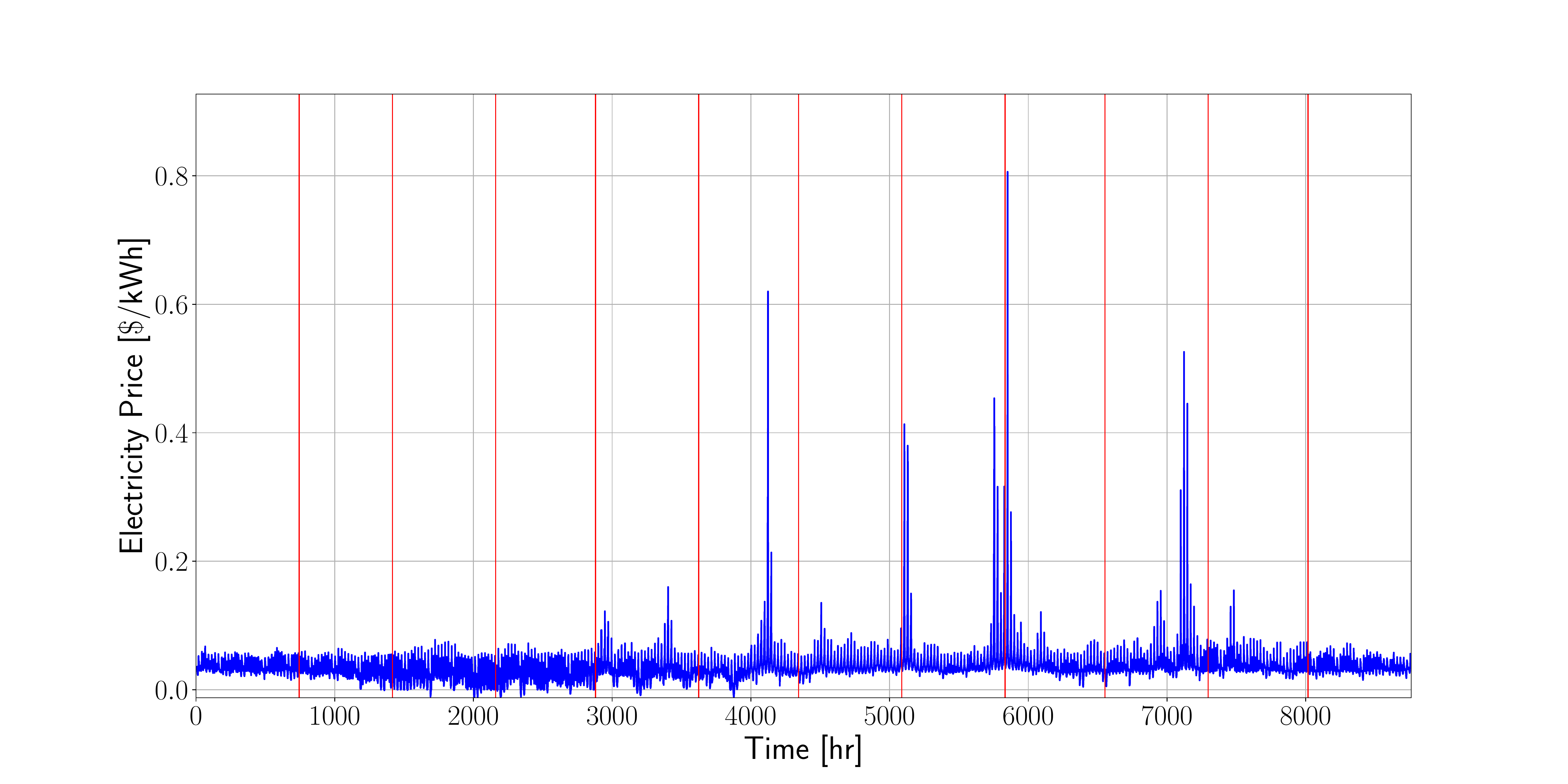}
\caption{Historical electricity price data. Red vertical lines denote end of each month.}
\label{fig:eprice}
\end{figure}

Figures \ref{fig:eload}-\ref{fig:eprice} show historical data for the campus electrical load, hot water load, chilled water load, and the electricity prices for the entire year. The vertical red lines in these figures represent monthly periods. In this work, we use these data as the available forecasts of the disturbances in the MPC formulation.

\FloatBarrier
\FloatBarrier
\FloatBarrier
\FloatBarrier

\subsection{Results} \label{subsec:results}

We now present computational results for the mixed-integer MPC problem \eqref{eq:detmpc}. We first use the DDIP scheme to solve the LP relaxation of this problem that represents the problem with zero minimum capacity of all units and show that the scheme obtains the optimal solution within an optimality gap of 0.01\%. We then present results for the mixed-integer MPC with different horizon lengths and show that the DDIP scheme provides high-quality solutions within optimality gaps of 0.1\%. The DDIP schemes are implemented in {\tt Julia} and leverage the algebraic modeling capabilities of {\tt JuMP} \cite{Dunning2017}. The stage-wise optimization subproblems in the forward and backward sweeps are solved using {\tt Gurobi 9.0}. We also use {\tt Gurobi} to solve the entire MPC problem (extensive form) and with this compare performance with the DDIP scheme. The algorithms were run on a 32-core machine with Ubuntu 14.04, Intel Xeon 2.30 GHz processors, and 188 GB RAM. 

\subsubsection{Zero Minimum Capacity (Continuous) MPC Problem}

\begin{figure}[!htp]
\centering
\includegraphics[width=0.7\textwidth]{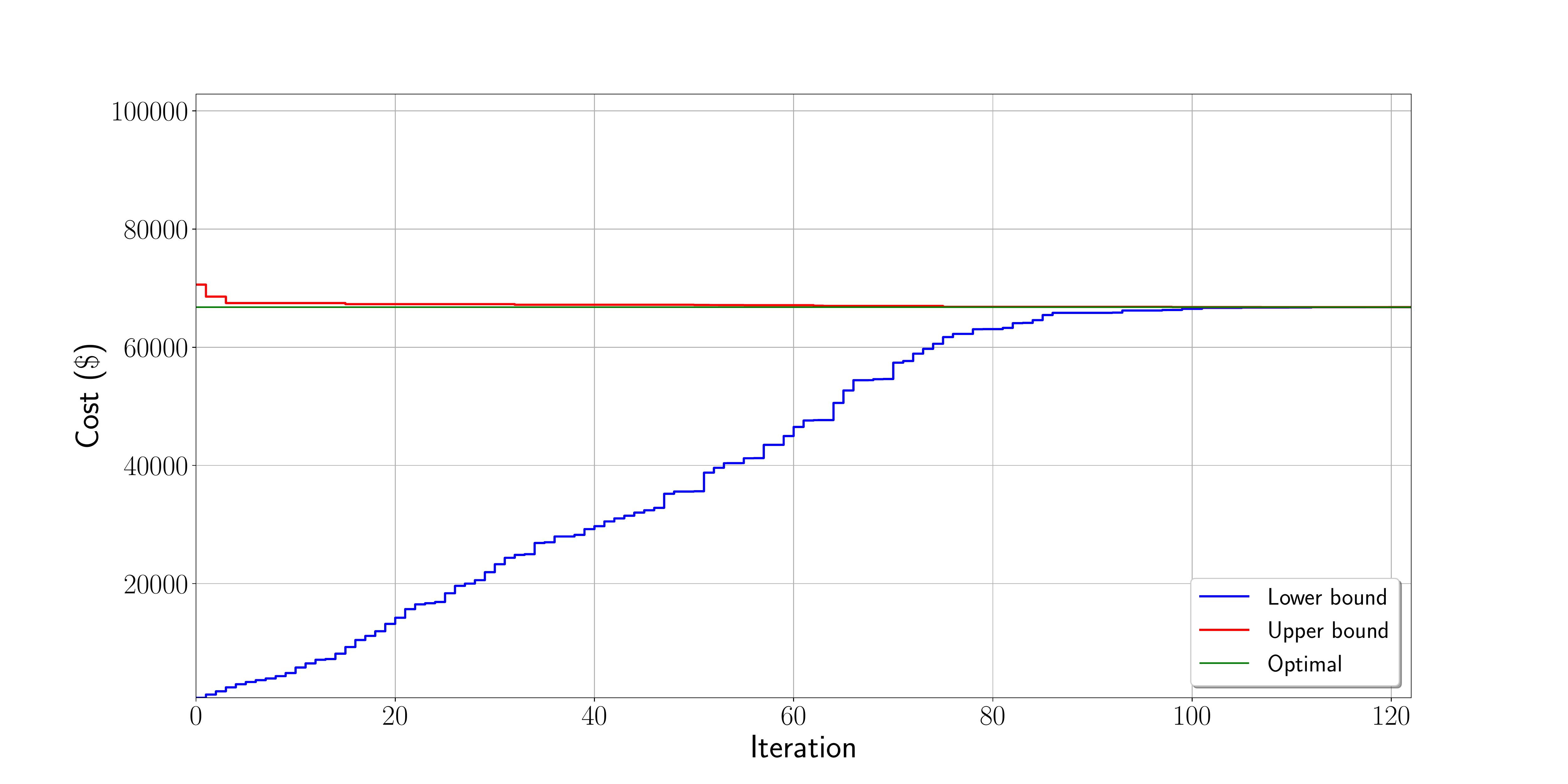}
\includegraphics[width=0.7\textwidth]{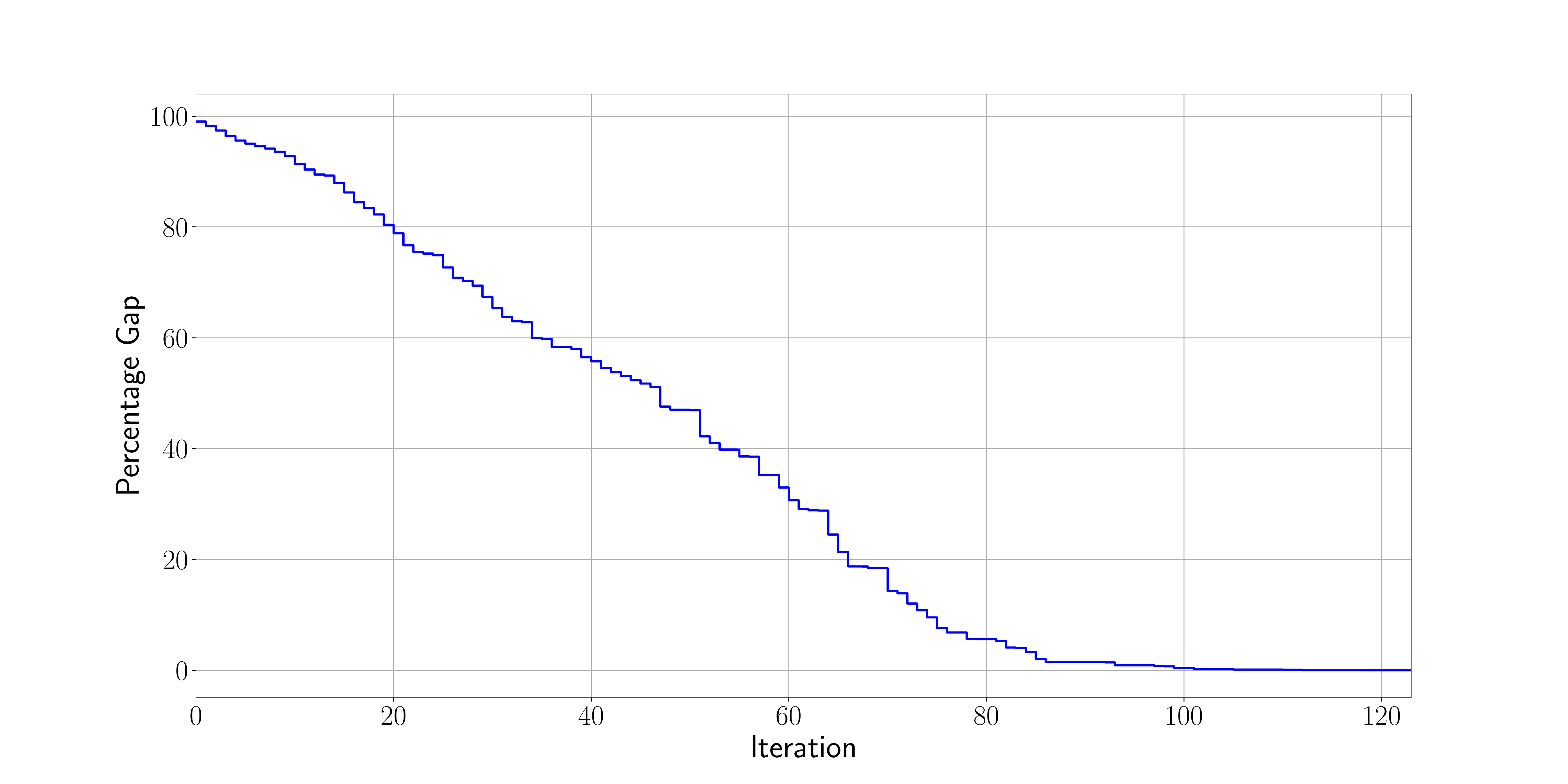}
\caption{\textbf{Top:} Evolution of upper and lower bounds over iterations for 168-hour horizon, continuous MPC; red curve represents the evolution of the upper bound and blue curve represents the evolution of the lower bound; green line represents the optimal objective value. \textbf{Bottom:} Evolution of the percentage optimality gap over iterations for 168-hour horizon, continuous MPC.}
\label{fig:gap_168_cont}
\end{figure}

\begin{figure}[!ht]
\centering
\includegraphics[width=0.7\textwidth]{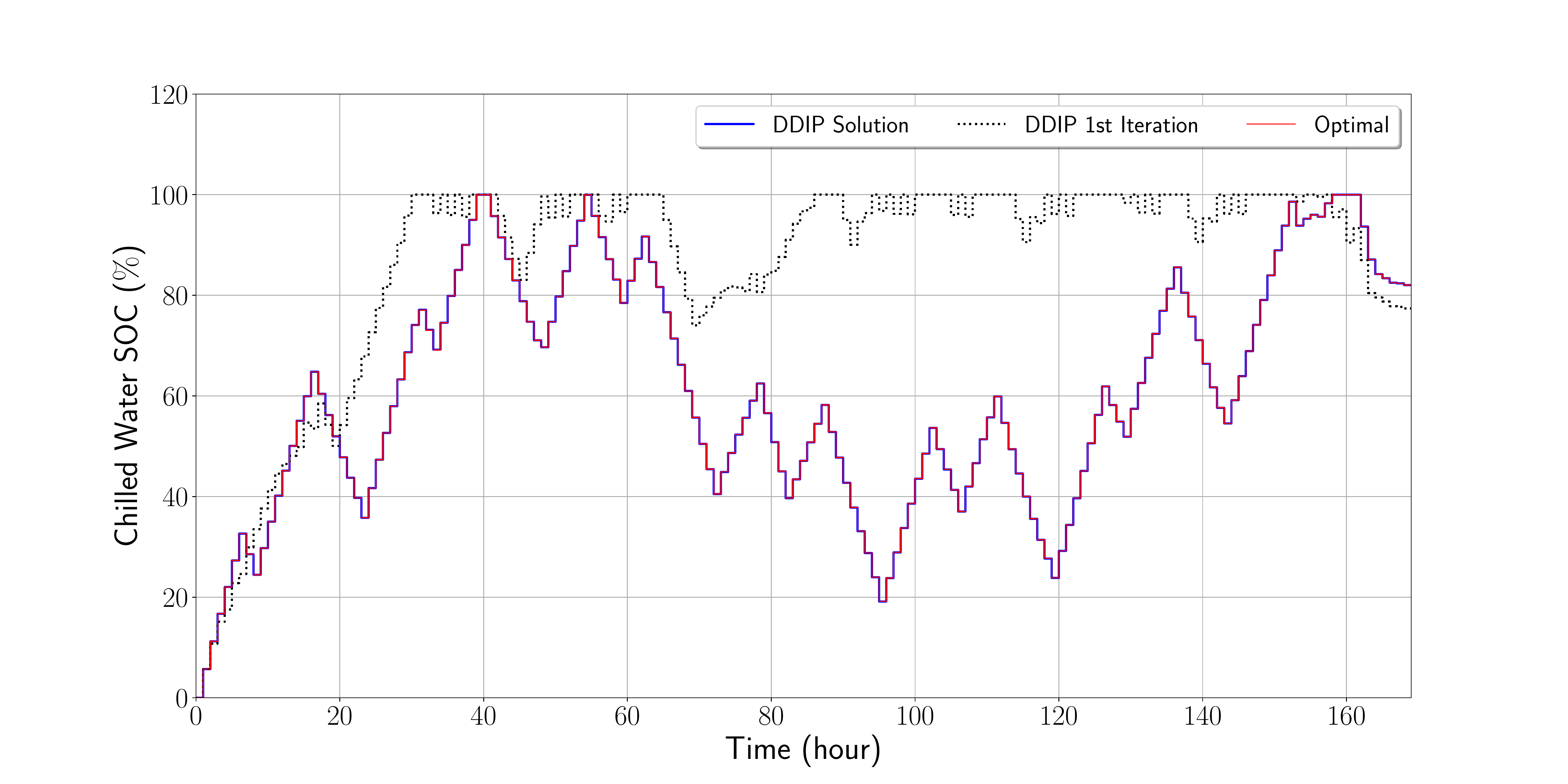}
\includegraphics[width=0.7\textwidth]{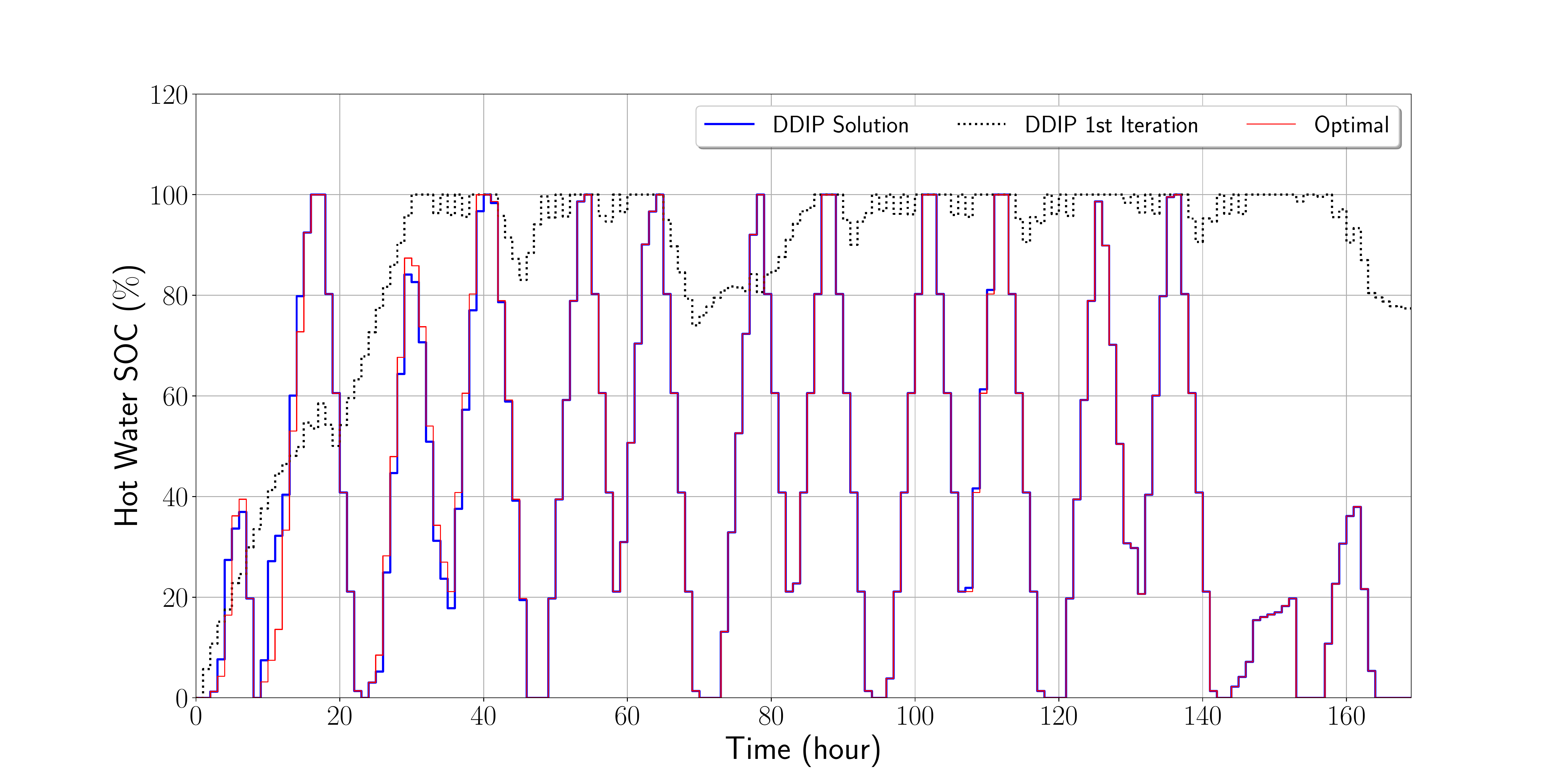}
\caption{Comparison of DDIP and optimal solution of the SOC evolution for chilled water tank (top) and hot water tank (bottom) for 168-hour horizon, continuous MPC. Blue curve represents the solution from the DDIP scheme; red curve represents the optimal solution; dotted-black curve represents the solution obtained at the first iteration of DDIP.}
\label{fig:soc_168_cont}
\end{figure}

We solve the MPC problem with zero minimum capacity of all units which is the LP relaxation of \eqref{eq:detmpc} defined over a time horizon with 168 timesteps (1 week with 1 hr resolution). The problem is an LP with 9,412 variables and 9,243 linear constraints. We implement the multiscale DDIP scheme with 84 stages and 2 internal timesteps within each stage. Figures \ref{fig:gap_168_cont}-\ref{fig:soc_168_cont} summarize the results. The top and bottom panels in Figure \ref{fig:gap_168_cont} show that the DDIP scheme obtains the optimal solution (within an optimality gap of 0.01\%) in 123 iterations. A total CPU time of 9.7s is required to achieve the stopping criteria (gap of 0.01\%).  We note that the extensive form of this LP problem is not computationally challenging and can be solved quite easily with the off-the-shelf solver \texttt{Gurobi} in 0.03s. As such, these results only seek to highlight that the DDIP scheme is consistent. The objective value obtained by DDIP solution (\$66,788.51) is close to the optimal objective (\$66,788.12). From Figures \ref{fig:gap_168_cont} we also see that the upper bound converges quickly to the optimal solution but the lower bound has slow convergence. We also see that the scheme consistently reduces the optimality gap (the relative difference between upper and lower bounds). 

\FloatBarrier

Figure \ref{fig:soc_168_cont} provides a comparison of the optimal solution and the final solution of the DDIP scheme for the SOC evolution of the chilled water and hot water storage tanks. It can be observed that the DDIP scheme obtains a close-to-optimal evolution of the states. The state evolution obtained at the first iteration of the DDIP scheme is also shown here as a dotted black curve. Here, we recall that the first iteration does not use cost-to-information and thus corresponds to a receding horizon scheme that neglects future information. It can be observed that the state evolution at the first DDIP iteration is far from the optimal evolution,  indicating that considering future information is important. The cutting plane information collected over the DDIP iterations close the optimality gap by using future information, bringing the state evolution closer to the optimal evolution. 

\subsubsection{Short-Horizon Mixed-Integer MPC Problem} 

We now use DDIP to solve a small mixed-integer MPC problem with a horizon comprising 168 timesteps. This is an MILP with 9,412 variables (including 3,360 binary variables) and 9,243 linear constraints. We implement the multiscale DDIP scheme using 84 stages (2 internal timesteps within each stage) and stop the algorithm when a 0.1\% gap between lower and upper bounds is reached. 

\begin{figure}[!htp]
\centering
\includegraphics[width=0.7\textwidth]{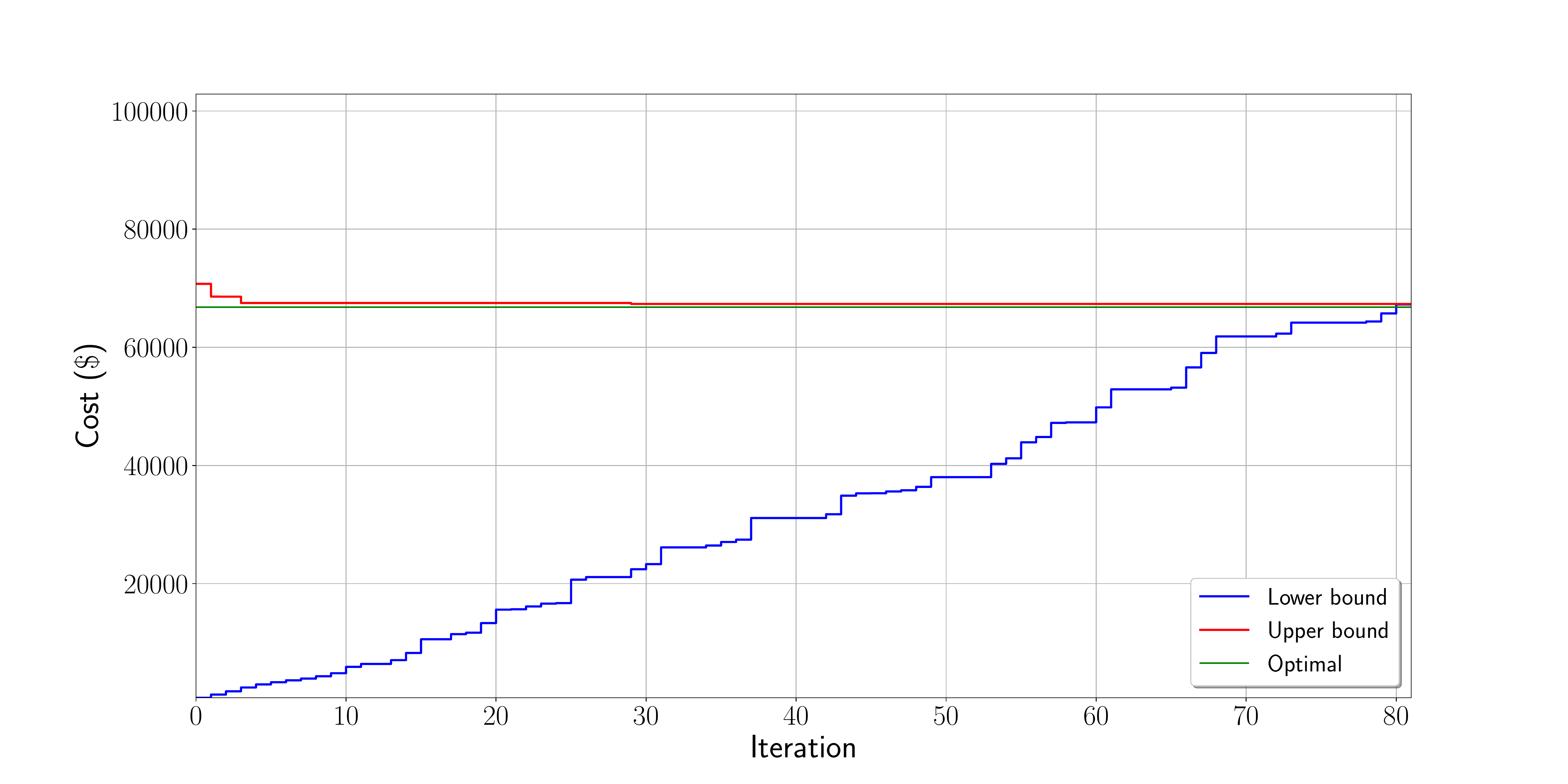}
\includegraphics[width=0.7\textwidth]{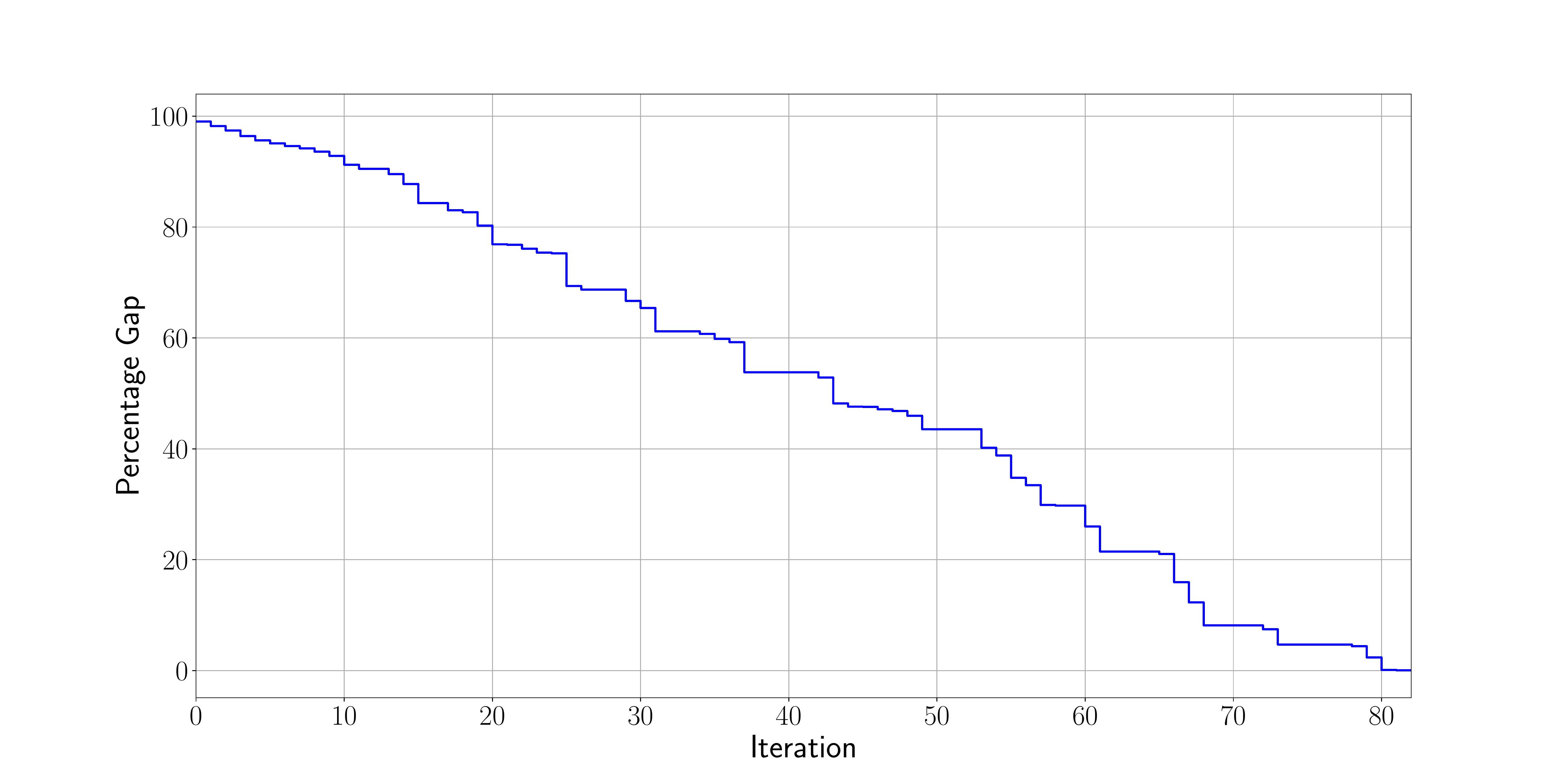}
\caption{\textbf{Top:} Evolution of DDIP upper and lower bounds for 168-hour horizon, mixed-integer MPC; red curve represents the evolution of the upper bound and blue curve represents the evolution of the lower bound; green line represents the optimal objective value. \textbf{Bottom:} Evolution of optimality gap for 168-hour horizon, mixed-integer MPC.}
\label{fig:gap_168_int}
\end{figure}

\begin{figure}[!htp]
\centering
\includegraphics[width=0.7\textwidth]{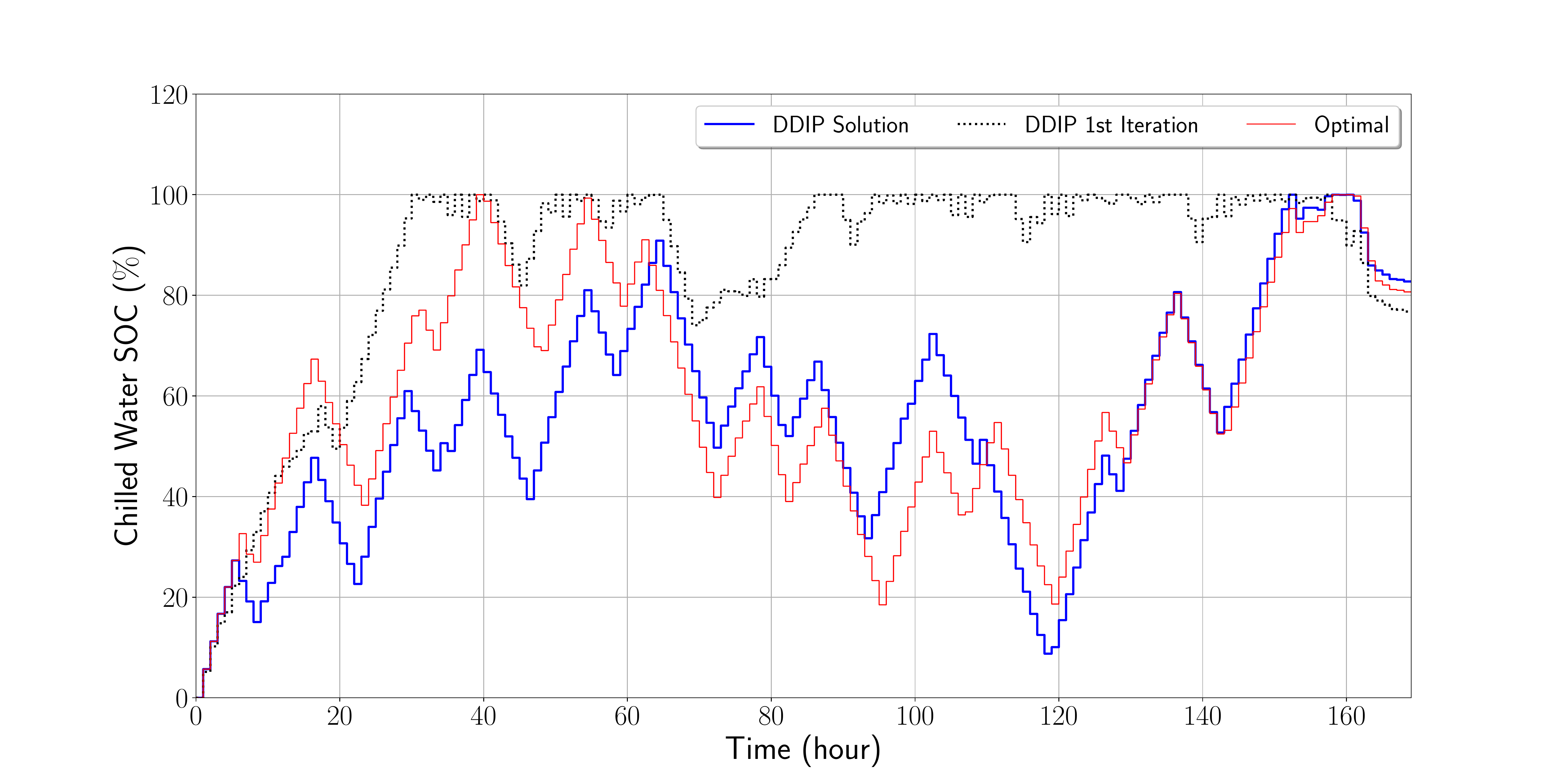}
\includegraphics[width=0.7\textwidth]{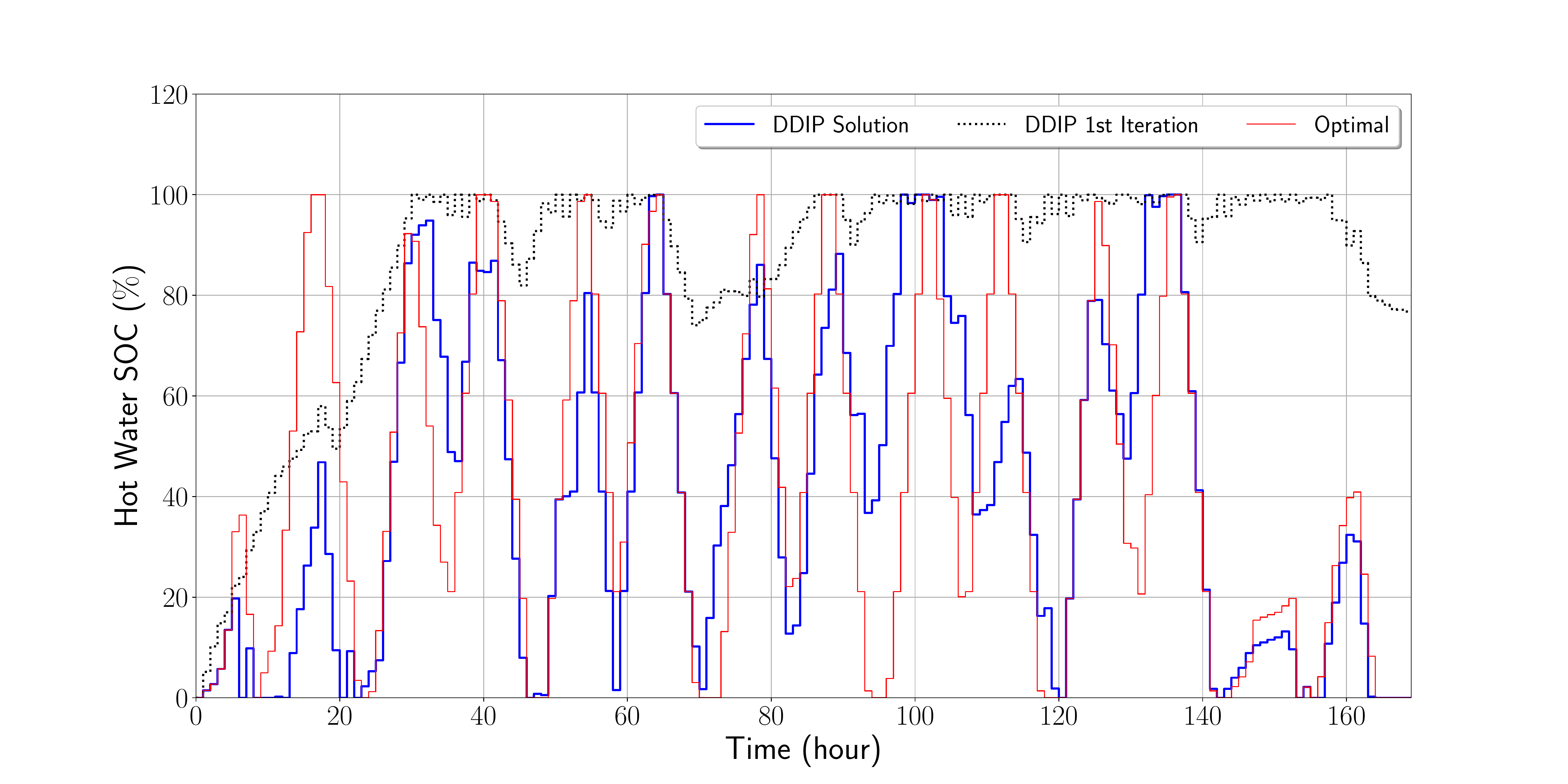}
\caption{Comparison of SOC evolution for chilled water tank (top) and hot water tank (bottom) for 168-hour horizon, mixed-integer MPC problem. Blue curve represents the solution from the DDIP scheme; red curve represents the optimal solution; dotted-black curve represents the solution obtained at the first iteration of DDIP.}
\label{fig:soc_168_int}
\end{figure}

Figures \ref{fig:gap_168_int}-\ref{fig:soc_168_int} summarize the results. The top and bottom panels in Figure \ref{fig:gap_168_int} show that the DDIP scheme reaches the desired gap in approximately 80 iterations and 24s of CPU time. The solution time of DDIP is ten times larger than the solution time of \texttt{Gurobi} (2s); this shows that temporal decomposition does not pay off in this small problem. This also highlights that off-the-shelf solvers are highly efficient. The objective value obtained by DDIP (\$66,340.64) deviates from the optimal objective (\$66,791.86), this gives an optimality gap 0.89\%. We also highlight that the DDIP scheme closes the optimality gap to approximately 2.67\% in just a couple of iterations, and it takes the remaining time to close the gap within the stopping criteria. This indicates that high-quality information about the future can be quickly obtained with a few backward sweeps.  

\begin{figure}[!htp]
\centering
\includegraphics[width=0.7\textwidth]{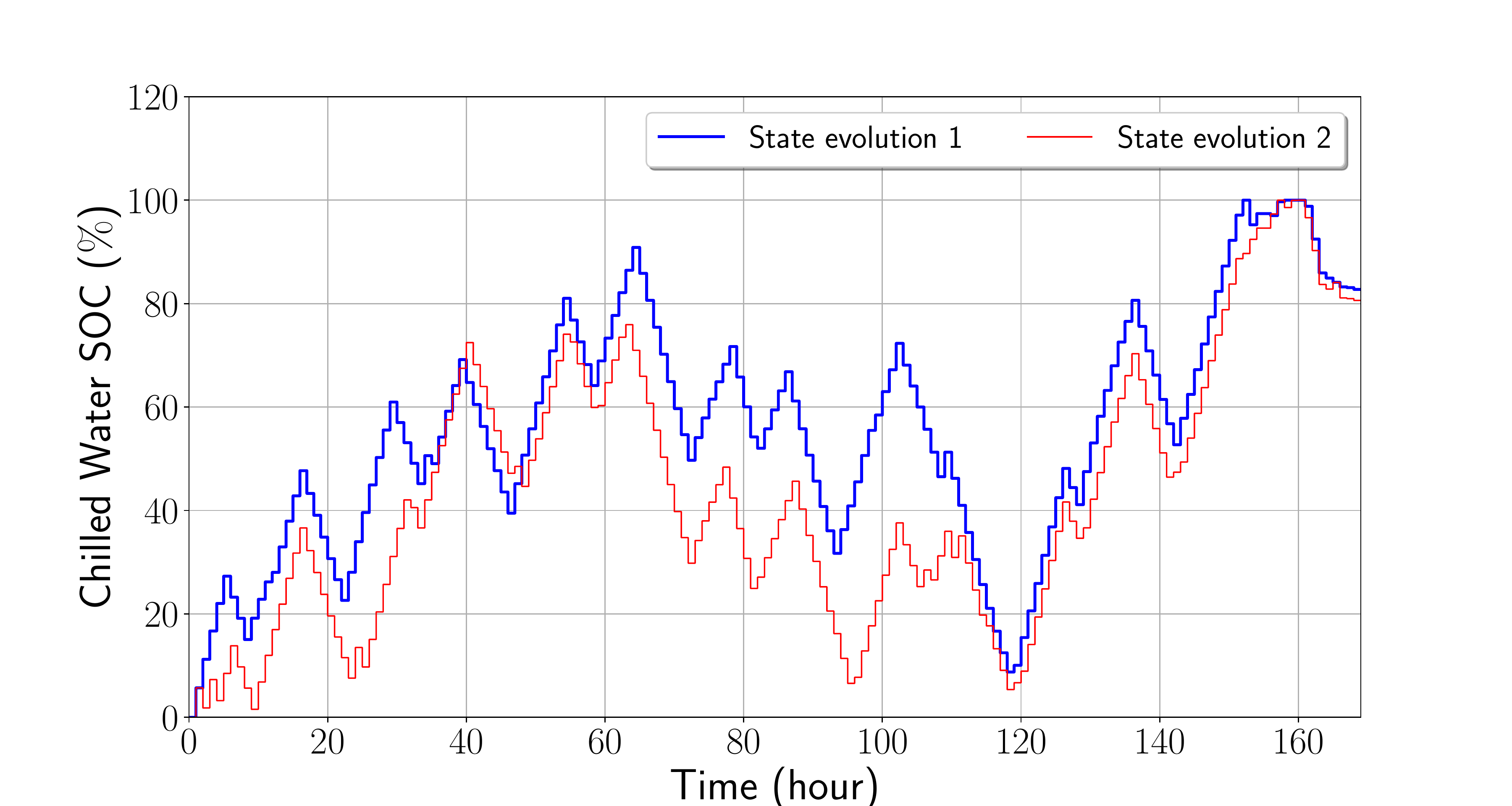}
\includegraphics[width=0.7\textwidth]{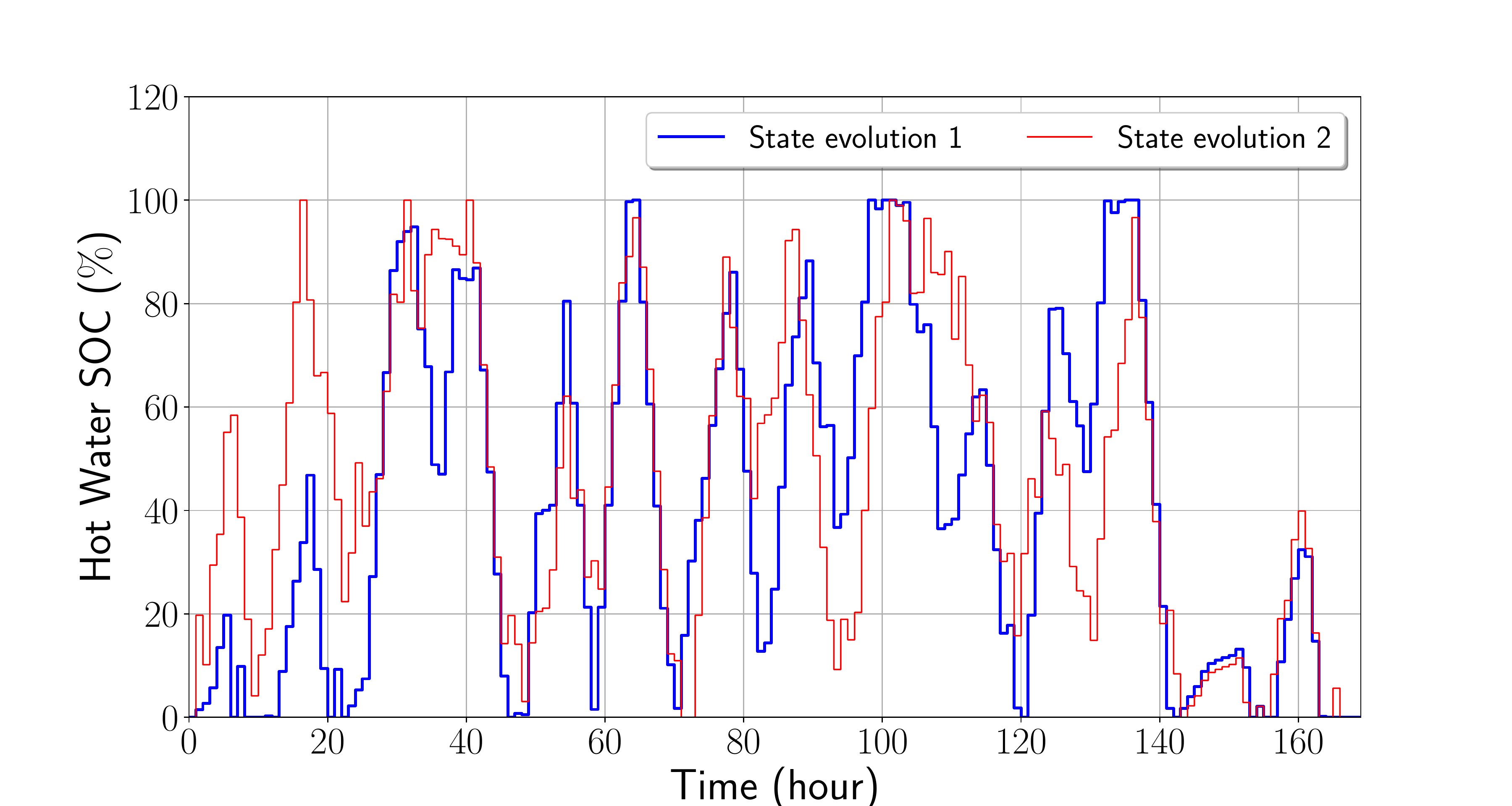}
\caption{Alternative solutions for the SOC evolution of chilled water tank (top) and hot water tank (bottom) achieving the same optimal objective value of \$67,340.64 for the 168-hour horizon, mixed-integer MPC problem. Both solutions are obtained with {\tt Gurobi}.}
\label{fig:soc_168_int_degen}
\end{figure}

Figure \ref{fig:soc_168_int} provides a comparison of the optimal solution and the final solution from the DDIP scheme for the SOC evolution of the chilled water and hot water tanks. It can be observed that the DDIP scheme obtains a similar to the evolution of the states of charge but larger deviations are obtained (compared to the LP case). Deviations from the optimal solution are attributed to the larger gap observed (which is small but trigger deviations in the optimal policy). The dotted black curve shows that the solution obtained at the first iteration of the DDIP scheme; we see that this policy is far from the optimal state evolution. Another potential reason for the observed deviation is that the mixed-integer MPC problem is degenerate (it contains multiple solutions with the same objective value). In Figure \ref{fig:soc_168_int_degen} we show that, {\em for the same optimal objective value} of \$67,340, there are at least two distinct SOC evolutions for the hot water and chilled water tanks. These solutions are obtained with Gurobi using no-good cuts. It is thus possible that the DDIP algorithm stops at a solution that has a different state evolution than the optimal solution obtained by \texttt{Gurobi} (although having close objective values). 

\FloatBarrier

\subsubsection{Long-Horizon Mixed-Integer MPC Problem}

To analyze the scalability of the DDIP scheme, we solve the mixed-integer MPC problem for the central HVAC plant over horizons spanning 2, 4, 5, 8, 10, and 20 weeks (all with 1-hour resolutions). The largest problem solved contains $20\times 168=3,360$ hourly timesteps. The results are summarized in Tables \ref{tab:first_three} and \ref{tab:gurobi_comparison}. For all problems, we implement the multiscale DDIP scheme with the same number stages (84) and vary the number of inner timesteps. 

\begin{figure}[!htp]
\centering
\includegraphics[width=0.7\textwidth]{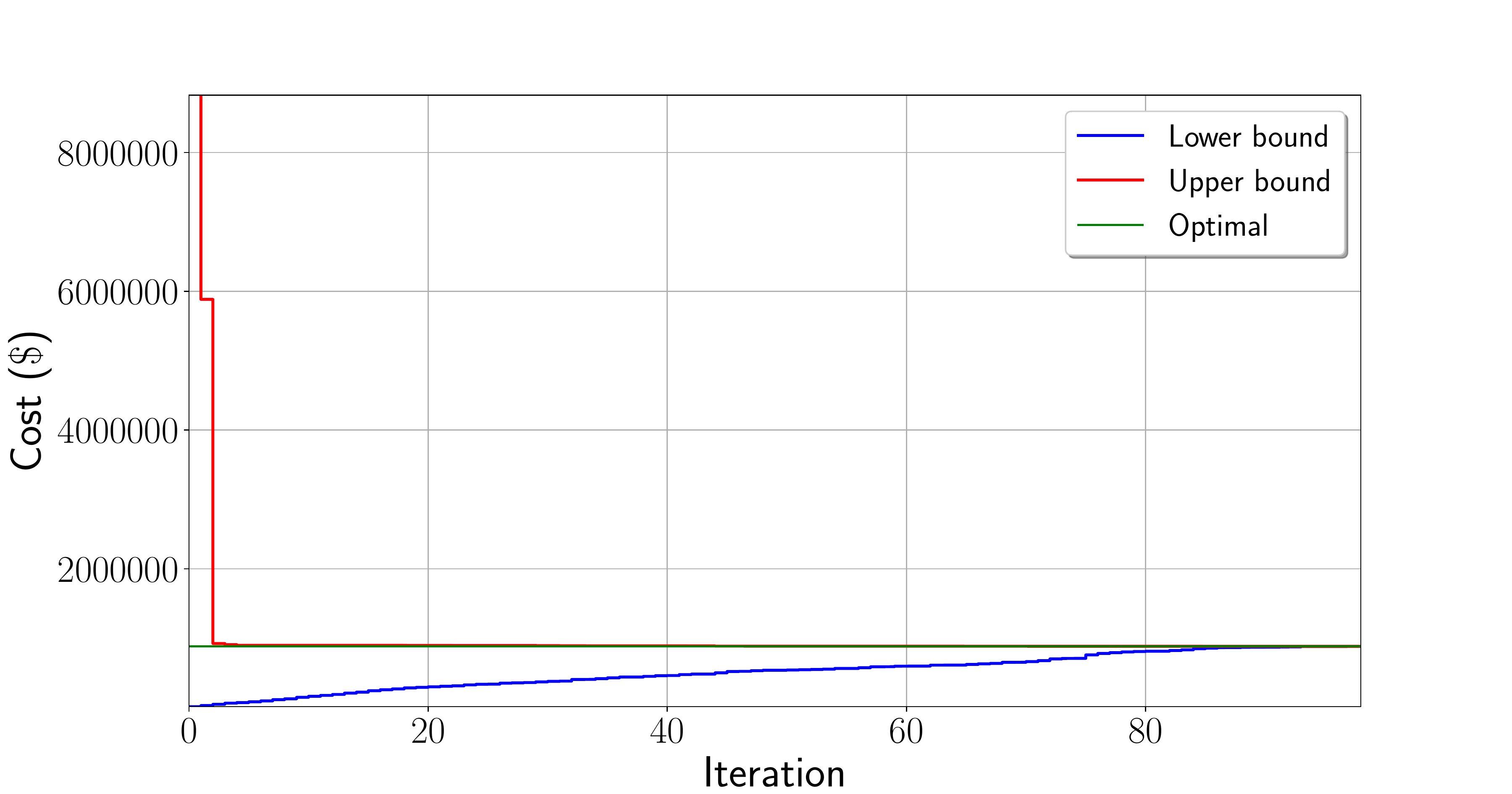}
\caption{Evolution of DDIP upper and lower bounds for 3360-hour horizon, mixed-integer MPC; showing in red curve represents the evolution of the upper bound and blue curve represents the evolution of the lower bound; green line represents the optimal objective value.}
\label{fig:gap_3360_int}
\end{figure}

In Table \ref{tab:first_three} we can see that the first iteration of the DDIP scheme delivers an optimality gap  (best upper bound of DDIP relative to the optimal solution obtained with \texttt{Gurobi}) of 5\% (except for the large problem). This indicates that, for the small problems, the value of future information is relatively small and thus the receding horizon policy provides a good approximation. For the large problem, however, we see that the first iteration has an optimality gap of over 4000\%. This is because the receding horizon policy induces infeasibility in satisfying the demands (which are heavily penalized). This indicates that, in this case, the value of future information is significant. In this problem, we also see that after 3 iterations, the DDIP scheme reduces the optimality gap to less than 5\%; this indicates that cutting planes are effective at removing the infeasibility. Figure \ref{fig:gap_3360_int} shows the evolution of the upper and lower bounds of the DDIP iterations for this problem. The huge reduction in the upper bound in the first 3 iterations can be clearly observed and shows that infeasibility is removed in the initial phase of the DDIP iterations.

\FloatBarrier

Table \ref{tab:gurobi_comparison} provides a comparison of the solutions obtained from \texttt{Gurobi} and the multiscale DDIP scheme along with the computational times. We observe that, for smaller-sized problems with horizons of up to 5 weeks, the DDIP scheme requires more time than \texttt{Gurobi}. However, DDIP outperforms {\tt Gurobi} for problems with 8 weeks and more. Specifically, we highlight that, for the 20-week long horizon problem, \texttt{Gurobi} requires a CPU time of approximately 7 hours (420 min), while the DDIP scheme solves the problem in 11 min within an optimality gap of 0.14\%. This represents a reduction in computational time of 97\%. We also observe that the number of iterations required to achieve the stopping criteria below 0.1\% gap (between upper and lower bounds) is consistently less than 100. As a result, the time required to achieve convergence scales close to linearly with the horizon length. For example, for the 2-week problem, DDIP requires half the solution time of the 4-week problem. The solution time for the 20-week problem is 27 times the solution time of 1-week problem. This highlights the desirable scalability properties of DDIP.  From Table \ref{tab:gurobi_comparison}  we also highlight that the largest problem solved contains 188,164 variables (67,200 of these are binary) and 184,803 constraints.  From this table, we can also see the fast explosion in the computational time of {\tt Gurobi} as we increase the horizon length. This is the result of an increasingly large combinatorial space.

\begin{table*}[!htb]\centering
\caption{Summary of solutions obtained by \texttt{Gurobi} and initial iterations of DDIP for the mixed-integer MPC problem.}
\resizebox{\columnwidth}{!}{%
\begin{tabular}{|c|c|c|c|c|c|}
\hline
\multirow{2}{*}{\textbf{Problem}} & {\textbf{\texttt{Gurobi}}} & \multicolumn{3}{c|}{\textbf{DDIP}}\\
\cline{2-5}
{} & {\begin{tabular}{@{}c@{}}{Optimal} \\{Objective (\$)}\end{tabular}} & {Iteration} &  {UB (\$)} & {\begin{tabular}{@{}c@{}}{\% Optimality Gap} \\ {$\mods{UB-\text{opt. obj.}}$} \end{tabular}}  \\
\hline
\multirow{4}{*}{ {\begin{tabular}{@{}c@{}}{Horizon: 168} \\{Internal timesteps: 2} \\ {Stages: 84}\end{tabular}} } & \multirow{4}{*}{66,792} & {1} & {70,742} & {5.91\%}  \\
\cline{3-5}
{} & {} & {2} & {68,576} & {2.67\%}  \\
\cline{3-5}
{} & {} & {3} & {68,575} & {2.67\%}  \\
\cline{3-5}
{} & {} & {Best} & {67,341} & {0.89\%}  \\
\hline
\multirow{4}{*}{ {\begin{tabular}{@{}c@{}}{Horizon: $2 \times 168$} \\{Internal timesteps: 4} \\ {Stages: 84}\end{tabular}} } & \multirow{4}{*}{120,335} & {1} & {127,848} & {6.24\%} \\
\cline{3-5}
{} & {} & {2} & {123,788} & {2.87\%} \\
\cline{3-5}
{} & {} & {3} & {123,271} & {2.44\%} \\
\cline{3-5}
{} & {} & {Best} & {121,350} & {0.84\%} \\
\hline
\multirow{4}{*}{ {\begin{tabular}{@{}c@{}}{Horizon: $4 \times 168$} \\{Internal timesteps: 8} \\ {Stages: 84}\end{tabular}} } & \multirow{4}{*}{256,062} & {1} & {269,252} & {5.15\%} \\
\cline{3-5}
{} & {} & {2} & {259,307} & {1.27\%} \\
\cline{3-5}
{} & {} & {3} & {259,307} & {1.27\%} \\
\cline{3-5}
{} & {} & {Best} & {256,899} & {0.33\%} \\
\hline
\multirow{4}{*}{ {\begin{tabular}{@{}c@{}}{Horizon: $5 \times 168$} \\{Internal timesteps: 10} \\ {Stages: 84}\end{tabular}} } & \multirow{4}{*}{302,842} & {1} & {315,963} & {4.33\%} \\
\cline{3-5}
{} & {} & {2} & {312,246} & {3.11\%} \\
\cline{3-5}
{} & {} & {3} & {308,412} & {1.84\%} \\
\cline{3-5}
{} & {} & {Best} & {304,167} & {0.44\%} \\
\hline
\multirow{4}{*}{ {\begin{tabular}{@{}c@{}}{Horizon: $8 \times 168$} \\{Internal timesteps: 16} \\ {Stages: 84}\end{tabular}} } & \multirow{4}{*}{414,919} & {1} & {434,402} & {4.70\%} \\
\cline{3-5}
{} & {} & {2} & {434,402} & {4.70\%} \\
\cline{3-5}
{} & {} & {3} & {429,337} & {3.47\%} \\
\cline{3-5}
{} & {} & {Best} & {416,265} & {0.32\%} \\
\hline
\multirow{4}{*}{ {\begin{tabular}{@{}c@{}}{Horizon: $10 \times 168$} \\{Internal timesteps: 20} \\ {Stages: 84}\end{tabular}} } & \multirow{4}{*}{484,720} & {1} & {507,866} & {4.78\%} \\
\cline{3-5}
{} & {} & {2} & {507,866} & {4.78\%} \\
\cline{3-5}
{} & {} & {3} & {499,579} & {3.07\%} \\
\cline{3-5}
{} & {} & {Best} & {485,823} & {0.23\%} \\
\hline
\multirow{4}{*}{ {\begin{tabular}{@{}c@{}}{Horizon: $20 \times 168$} \\{Internal timesteps: 40} \\ {Stages: 84}\end{tabular}} } & \multirow{4}{*}{881,432} & {1} & {36,598,168} & {4052.13\%}  \\
\cline{3-5}
 {} & {} & {2} &  {5,883,037} & {567.44\%} \\
 \cline{3-5}
 {} & {} &  {3} &  {923,617} & {4.78\%} \\
 \cline{3-5}
 {} & {} &  {Best}  & {882,668} & {0.14\%} \\
\hline
\end{tabular} \label{tab:first_three}
}
\end{table*}

\begin{table*}[!htb]\centering
\caption{Comparison of solutions and computational times of \texttt{Gurobi} and DDIP scheme.}
\resizebox{\columnwidth}{!}{%
\begin{tabular}{|c|c|c|c|c|c|c|c|c|c|c|c|}
\hline
\multirow{2}{*}{\textbf{Horizon}} & \multicolumn{3}{c|}{\textbf{Problem Size}} & \multicolumn{2}{c|}{\textbf{\texttt{Gurobi}}}  & \multicolumn{6}{c|}{\textbf{DDIP}}\\
\cline{2-12}
{\textbf{(Hours)}} & {\begin{tabular}{@{}c@{}}{Total} \\ {Variables}\end{tabular}} & {\begin{tabular}{@{}c@{}}{Binary} \\{Variables}\end{tabular}} & {Constraints} & {\begin{tabular}{@{}c@{}}{Optimal} \\{Objective (\$)}\end{tabular}} & {\begin{tabular}{@{}c@{}}{Solution} \\{Time}\end{tabular}}  & {\begin{tabular}{@{}c@{}}{Best} \\{UB (\$)}\end{tabular}} & {LB (\$)} &  {\begin{tabular}{@{}c@{}}{\% Gap} \\{$\mods{UB-LB}$}\end{tabular}} & {\begin{tabular}{@{}c@{}}{\% Optimality Gap} \\{$\mods{UB-\text{opt. obj.}}$}\end{tabular}} &  {Iterations} &  {\begin{tabular}{@{}c@{}}{Solution} \\{Time}\end{tabular}}  \\
\hline
{168} & {9,412} & {3,360} & {9,243} & {66,791.86} & {2 s} & {67,340.64} & {67,302.76} & {0.056\%} & {0.89\%} & {82} & {24 s} \\
\hline
{$2 \times 168$ } & {18,820} & {6,720} & {18,483} & {120,335.40} & {6 s} & {121,244.65} & {121,349.54} & {0.086\%} & {0.842\%} & {86} & {37 s}\\
\hline
{$4 \times 168$ } & {37,636} & {13,440} & {36,963} & {256,061.58} & {11 s} & {256,898.63} & {256,696.45} & {0.079\%} & {0.327\%} & {87} & {59 s}\\
\hline
{$5 \times 168$ } & {47,044} & {16,800} & {46,203} & {302,842.43} & {23 s} & {304,166.82} & {304,065.33} & {0.033\%} & {0.437\%} & {89} & {80 s}\\
\hline
{$8 \times 168$ } & {75,268} & {26,880} & {73,923} & {414,918.85} & {338 s} & {416,265.25} & {416,139.85} & {0.030\%} & {0.324\%} & {92} & {138 s}\\
\hline
{$10 \times 168$ } & {94,084} & {33,600} & {92,403} & {484,720.26} & {807 s} & {485,823.25} & {485,350.15} & {0.097\%} & {0.227\%} & {92} & {195 s}\\
\hline
{$20 \times 168$ } & {188,164} & {67,200} & {184,803} & {881,431.81} & {25,061 s} & {882,668.36} & {882,641.60} & {0.003\%} & {0.14\%} & {99} & {662 s} \\
\hline
\end{tabular} \label{tab:gurobi_comparison}
}
\end{table*}

\FloatBarrier

\begin{figure}[!htp]
\centering
\includegraphics[width=0.7\textwidth]{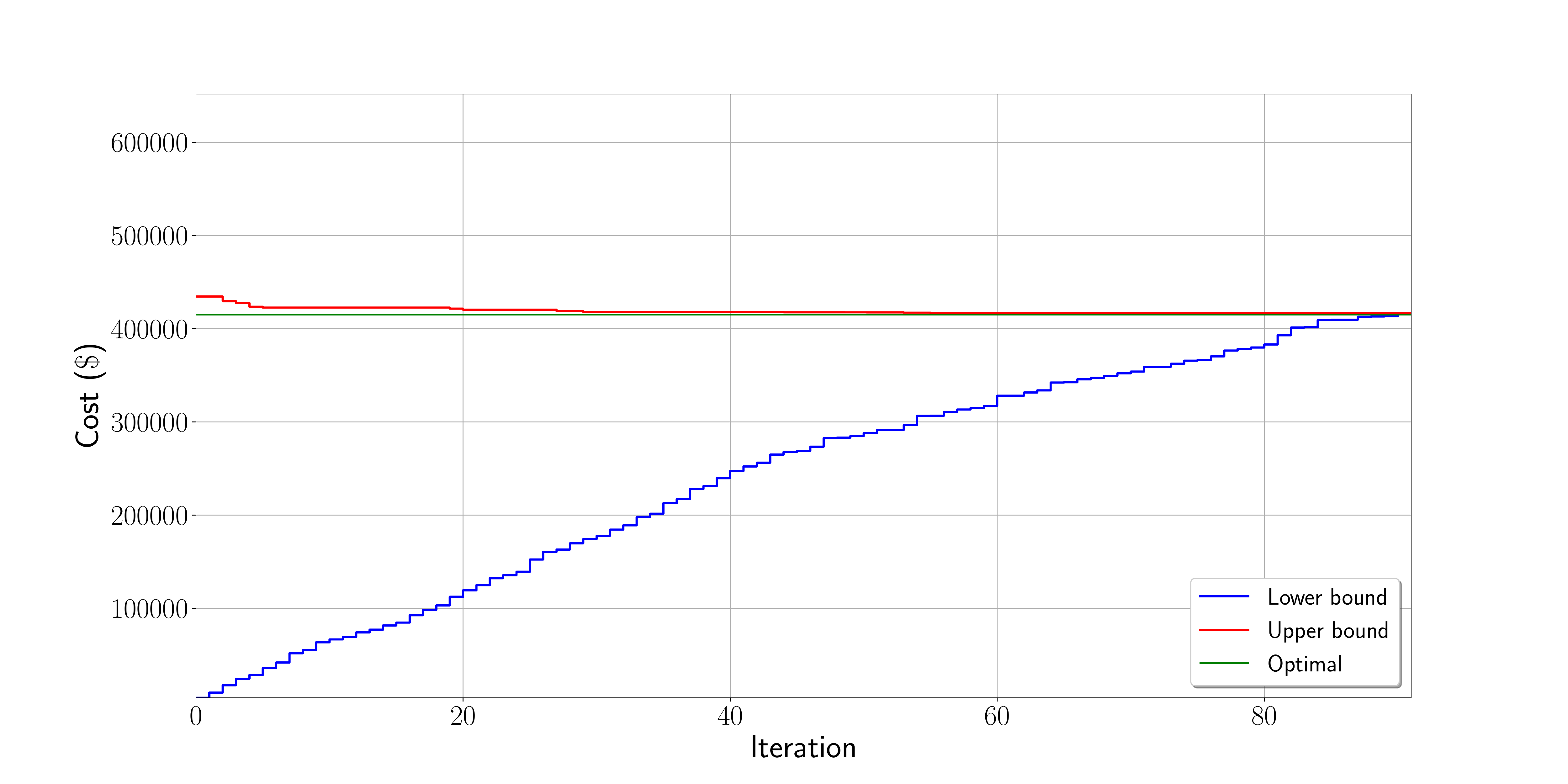}
\includegraphics[width=0.7\textwidth]{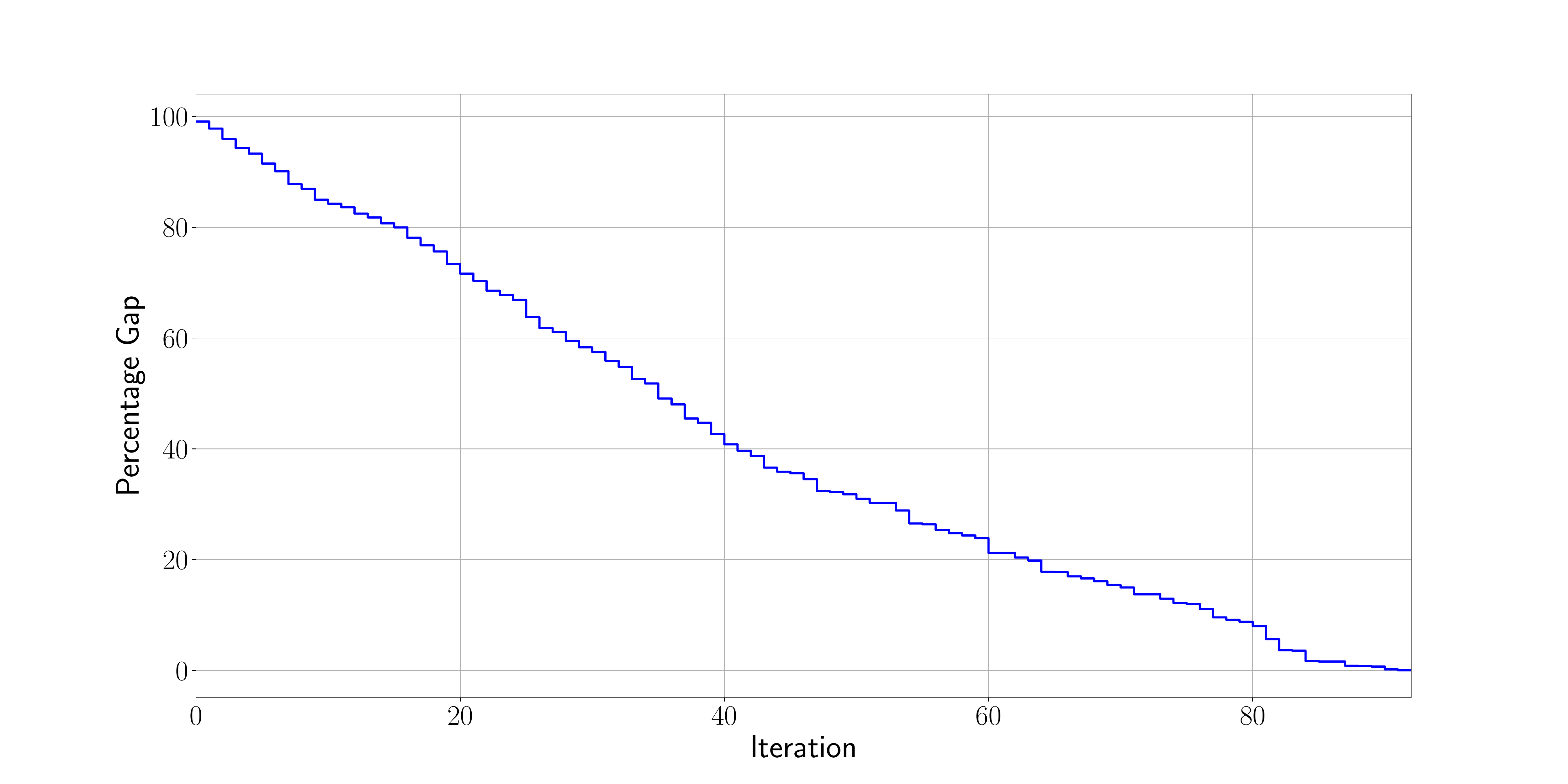}
\caption{\textbf{Top:} Evolution of upper and lower bounds over iterations for 1344-hour horizon, mixed-integer MPC; red curve represents the evolution of the upper bound and blue curve represents the evolution of the lower bound; green line represents the optimal objective value. \textbf{Bottom:} Evolution of the percentage optimality gap over iterations for 1344-hour horizon, mixed-integer MPC.}
\label{fig:gap_1344_int}
\end{figure}

\begin{figure}[!htp]
\centering
\includegraphics[width=0.7\textwidth]{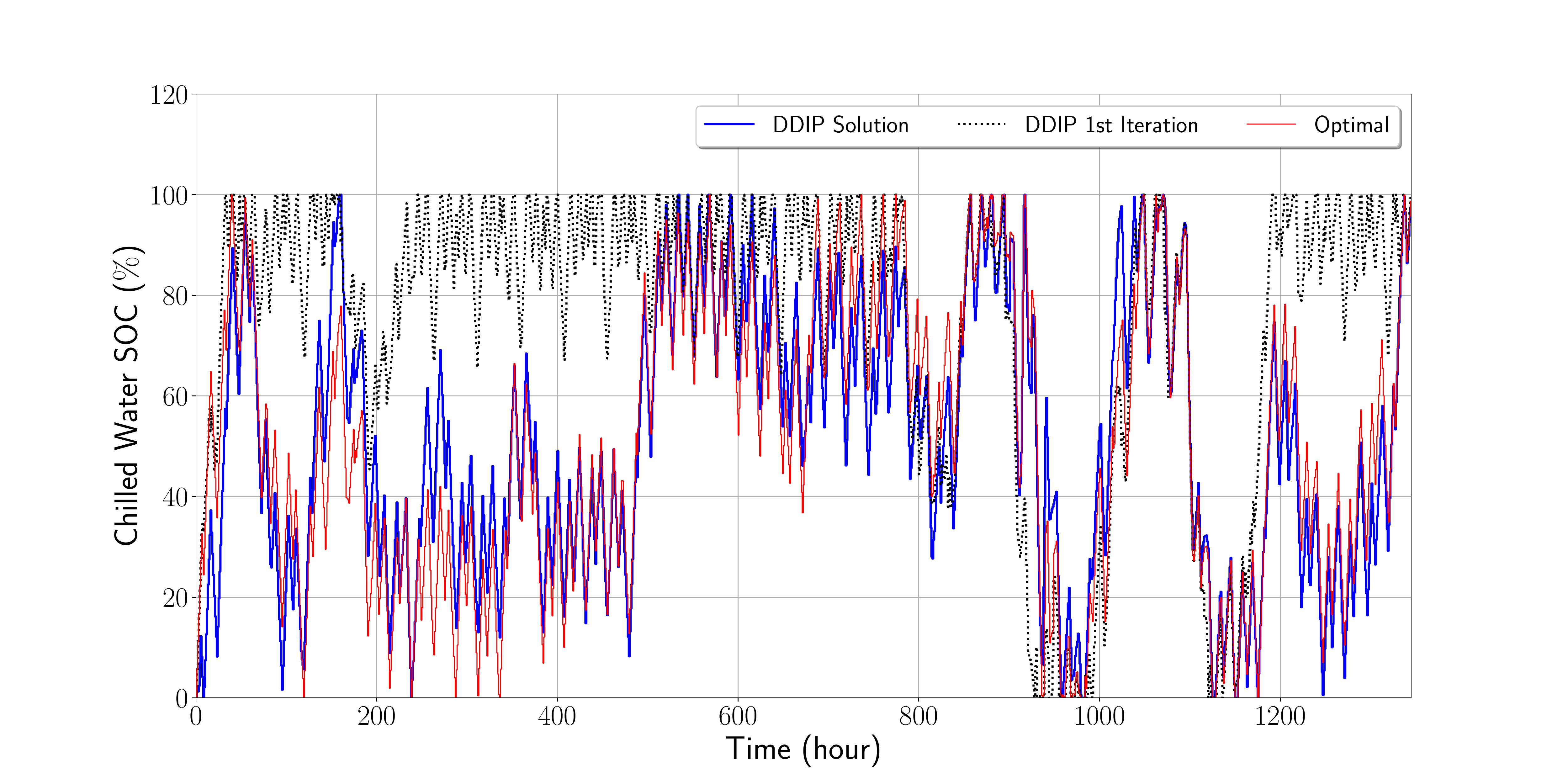}
\includegraphics[width=0.7\textwidth]{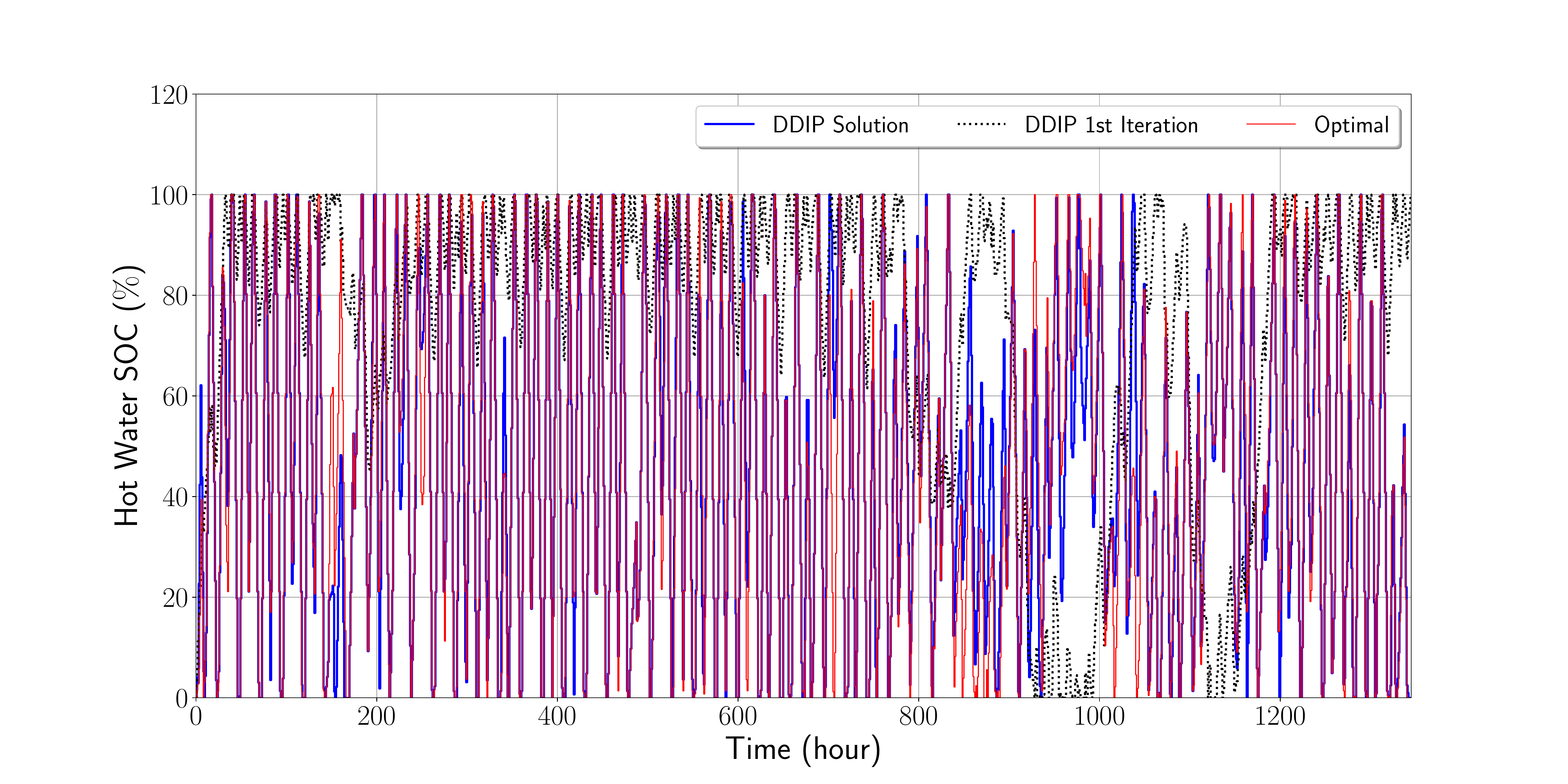}
\caption{Comparison of DDIP and optimal solution of the SOC evolution for chilled water tank (top) and hot water tank (bottom) for 840-hour horizon, mixed-integer MPC. Blue curve represents the solution from the DDIP scheme; red curve represents the optimal solution; dotted-black curve represents the solution obtained at the first iteration of DDIP.}
\label{fig:soc_1344_int}
\end{figure}

We further analyzed the performance of the multiscale DDIP scheme for a large mixed-integer MPC problem with 8-week (1344 timesteps). The full MILP problem contains 75,268 variables (26,880 binary variables) and 73,923 linear constraints. The top and bottom panels in Figure \ref{fig:gap_1344_int}-\ref{fig:soc_1344_int} present the results for this problem. From the top and bottom panels in Figure \ref{fig:gap_1344_int}, we see that the DDIP scheme achieves an optimality gap within the stopping criteria after 90 iterations. It requires a CPU time of 138s to achieve the stopping criteria of below 0.1\% gap, while \texttt{Gurobi} requires more than 300s to obtain the optimal solution. Thus the solution time by the multiscale DDIP scheme for the 8-week problem is less than 8 times that of the 1-week problem, which exhibits the scalability of the multiscale DDIP scheme. The final objective value obtained by the DDIP scheme is within a gap of 0.324\% from the optimal objective value. Since the multiscale DDIP scheme needs to solve small single-stage problems in every iteration, the algorithm easily scales to longer horizons.  Figure \ref{fig:soc_1344_int} compares an optimal solution obtained with {\tt Gurobi} and the final solution of the DDIP scheme for the SOC evolution of the chilled water and hot water tanks. It can be observed that the solution of the DDIP scheme follows the overall features of the optimal solution (which some high-frequency deviations). Here we also compare against the evolution after the first DDIP iteration; from here we see that the policy of the receding horizon scheme does not follow the general features of the policy. This, again, highlights the value of future information. 

\section{Conclusions and Future Work}
We proposed a dual dynamic integer programming (DDIP) scheme for solving multiscale mixed-integer MPC problems. The scheme uses a cutting-plane information collected over forward and backward sweeps along the time horizon to adaptively construct and refine cost-to-go functions. We show that the approach can handle general MPC formulations and that it can dramatically outperform state-of-the-art solvers. Future work will focus on the implementation of the proposed scheme for stochastic mixed-integer MPC applications with uncertainties evolving over long time horizons and/or fine time resolutions. We are also interested in exploring more general schemes that use parallelization and that use adaptive time discretizations. 

\section{Acknowledgments}

We acknowledge funding from the National Science Foundation under award NSF-EECS-1609183. VZ declares a financial interest in Johnson Controls International, a for-profit company that develops control technologies. We dedicate this paper to Sebastian Engell, a pioneer in the area of process systems engineering.

\appendix
\section*{APPENDIX}
\section{Nomenclature} \label{sec:nomencleture}
\noindent \textbf{Sets and indices:}
\begin{itemize}
\item[\textbullet] $\mathcal{T} :=  \{0,1,\dots,N\}$: Prediction horizon time set for the state evolution, where $N$ is the prediction horizon length.
\item[\textbullet] $\bar{\mathcal{T}} :=  \{0,1,\dots,N-1\}$: Prediction horizon time set for the control trajectory, where $N$ is the prediction horizon length.
\item[\textbullet] $t$: Time instant index.
\end{itemize}

\noindent \textbf{Model Parameters and Data:}
\begin{itemize}
\item[\textbullet] $\pi^{e}_t \in \mathbb{R}_+$: Electricity price [\$/kWh] over the time interval $[t,(t+1)]$.
\item[\textbullet] $\pi^{w}_t \in \mathbb{R}_+$: Price of water [\$/gal] over the time interval $[t,(t+1)]$.
\item[\textbullet] $\pi^{ng}_t \in \mathbb{R}_+$: Price of natural gas [\$/kWh] over the time interval $[t,(t+1)]$.
\item[\textbullet] $\alpha^{e}_{cs} \in \mathbb{R}_+$: kW of electricity used by each chiller unit per kW chilled water produced [-].
\item[\textbullet] $\alpha^{e}_{hrc} \in \mathbb{R}_+$: kW of electricity used by each HR chiller unit per kW chilled water produced [-].
\item[\textbullet] $\alpha^{e}_{hwg} \in \mathbb{R}_+$: kW of electricity used by each hot water generator per kW hot water produced [-].
\item[\textbullet] $\alpha^{e}_{ct} \in \mathbb{R}_+$: kW of electricity used by each cooling tower per kW condenser water input [-].
\item[\textbullet] $\alpha^{w}_{ct} \in \mathbb{R}_+$: Gallons of water used by each cooling tower per kW condenser water input [-].
\item[\textbullet] $\alpha^{ng}_{hwg} \in \mathbb{R}_+$: kW of natural gas used by each hot water generator per kW hot water produced [-].
\item[\textbullet] $\alpha^{cond}_{cs} \in \mathbb{R}_+$: kW of condenser water produced by each chiller unit per kW chilled water produced [-].
\item[\textbullet] $\alpha^{h}_{hrc} \in \mathbb{R}_+$: kW of hot water produced by each HR chiller unit per kW chilled water produced [-].
\item[\textbullet] $\rho^{cw} \in \mathbb{R}_+$: Penalty for unmet chilled water load [\$/kWh].
\item[\textbullet] $\rho^{hw} \in \mathbb{R}_+$: Penalty for unmet hot water load [\$/kWh].
\item[\textbullet] $\beta_{cs} \in [0,1] = 0.5$: Minimum operating load (fraction of maximum capacity) of each chiller unit when ON.
\item[\textbullet] $\beta_{hrc} \in [0,1] = 0.8$: Minimum operating load (fraction of maximum capacity) of each HR chiller unit when ON.
\item[\textbullet] $\beta_{hwg} \in [0,1] = 0.5$: Minimum operating load (fraction of maximum capacity) of each hot water generator when ON.
\item[\textbullet] $\beta_{ct} \in [0,1] = 0.5$: Minimum operating load (fraction of maximum capacity) of each cooling tower when ON.
\item[\textbullet] $\beta_{cs} \in [0,1] = 0$: Minimum operating load (fraction of maximum capacity) of dump heat exchanger when ON.
\item[\textbullet] $L^{e}_{t} \in \mathbb{R}_+$: Electrical load of campus [kW] over the time interval $[t,(t+1)]$.
\item[\textbullet] $L^{cw}_{t} \in \mathbb{R}_+$: Chilled water load [kW] over the time interval $[t,(t+1)]$.
\item[\textbullet] $L^{hw}_{t} \in \mathbb{R}_+$: Hot water load [kW] over the time interval $[t,(t+1)]$.
\item[\textbullet] $\overline{E}_{cw}\in \mathbb{R}_+$: Energy storage capacity of chilled water energy storage [kWh].
\item[\textbullet] $\overline{E}_{hw}\in \mathbb{R}_+$: Energy storage capacity of hot water energy storage [kWh].
\item[\textbullet] $\overline{p}_{cs}\in \mathbb{R}_+$: Maximum load of each chiller unit [kW].
\item[\textbullet] $\overline{p}_{hrc}\in \mathbb{R}_+$: Maximum load of each heat recovery (HR) chiller unit [kW].
\item[\textbullet] $\overline{p}_{hwg}\in \mathbb{R}_+$: Maximum load of each hot water generator [kW].
\item[\textbullet] $\overline{p}_{ct}\in \mathbb{R}_+$: Maximum load of each cooling tower [kW].
\item[\textbullet] $\overline{P}_{cw}\in \mathbb{R}_+$: Maximum discharging rate of chilled water energy storage [kW].
\item[\textbullet] $\overline{P}_{hw}\in \mathbb{R}_+$: Maximum discharging rate of hot water energy storage [kW].
\end{itemize}

\noindent \textbf{Controls:}
\begin{itemize}
\item[\textbullet] $Y_{n,cs,t} \in \{0,1\}$: Binary variable for $n^{th}$ chiller unit: ON if 1, OFF if 0 over the time interval $[t,(t+1)]$. 
\item[\textbullet] $Y_{n,hrc,t} \in \{0,1\}$: Binary variable for $n^{th}$ heat recovery (HR) chiller unit: ON if 1, OFF if 0 over the time interval $[t,(t+1)]$. 
\item[\textbullet] $Y_{n,hwg,t} \in \{0,1\}$: Binary variable for $n^{th}$ hot water generator: ON if 1, OFF if 0 over the time interval $[t,(t+1)]$. 
\item[\textbullet] $Y_{n,ct,t} \in \{0,1\}$: Binary variable for $n^{th}$ cooling tower: ON if 1, OFF if 0 over the time interval $[t,(t+1)]$. 
\item[\textbullet] $Y_{n,hx,t} \in \{0,1\}$: Binary variable for $n^{th}$ dump heat exchanger (HX) unit: ON if 1, OFF if 0 over the time interval $[t,(t+1)]$. 
\item[\textbullet] $p_{n,cs,t} \in \mathbb{R}_+$: Amount of chilled water produced by $n^{th}$ chiller unit [kW] over the time interval $[t,(t+1)]$. 
\item[\textbullet] $p_{n,hrc,t} \in \mathbb{R}_+$: Amount of chilled water produced by $n^{th}$ heat recovery (HR) chiller unit [kW] over the time interval $[t,(t+1)]$. 
\item[\textbullet] $p_{n,hwg,t} \in \mathbb{R}_+$: Amount of hot water produced by $n^{th}$ hot water generator [kW] over the time interval $[t,(t+1)]$. 
\item[\textbullet] $p_{n,ct,t} \in \mathbb{R}_+$: Amount of condenser water input to the $n^{th}$ cooling tower [kW] over the time interval $[t,(t+1)]$.
\item[\textbullet] $p_{n,hx,t} \in \mathbb{R}_+$: Amount of hot water input to the $n^{th}$ dump heat exchanger (HX) [kW] over the time interval $[t,(t+1)]$. 
\item[\textbullet] $P_{cs,t} \in \mathbb{R}_+$: Total amount of chilled water produced by chiller subplant (4 chiller units) [kW] over the time interval $[t,(t+1)]$. 
\item[\textbullet] $P_{hrc,t} \in \mathbb{R}_+$: Total amount of chilled water produced by heat recovery (HR) chiller subplant (3 HR chiller units) [kW] over the time interval $[t,(t+1)]$. 
\item[\textbullet] $P_{hwg,t} \in \mathbb{R}_+$: Total amount of hot water produced by 3 hot water generators [kW] over the time interval $[t,(t+1)]$. 
\item[\textbullet] $P_{hx,t} \in \mathbb{R}_+$: Total amount of hot water input to the dump heat exchanger (HX) [kW] over the time interval $[t,(t+1)]$. 
\item[\textbullet] $P_{cw,t} \in \mathbb{R}$: Net charge/discharge rate [kW] of the chilled water energy storage over the time interval $[t,(t+1)]$. If $P_{cw,t}>0$, the chilled water is being discharged and if $P_{cw,t}<0$ the chilled water is being charged.
\item[\textbullet] $P_{hw,t} \in \mathbb{R}$: Net charge/discharge rate [kW] of the hot water energy storage over the time interval $[t,(t+1)]$. If $P_{hw,t}>0$, the hot water is being discharged and if $P_{hw,t}<0$ the hot water is being charged.
\item[\textbullet] $r^{e}_{t}\in \mathbb{R}$: Residual electrical load [kW] over the time interval $[t,(t+1)]$.
\item[\textbullet] $r^{w}_{t}\in \mathbb{R}$: Residual water demand [gal/h] over the time interval $[t,(t+1)]$.
\item[\textbullet] $r^{ng}_{t}\in \mathbb{R}$: Residual natural gas demand [kW] over the time interval $[t,(t+1)]$.
\item[\textbullet] $U_{cw,t}\in \mathbb{R}$: Slack variable for unmet chilled water load [kW] over the time interval $[t,(t+1)]$.
\item[\textbullet] $U_{hw,t}\in \mathbb{R}$: Slack variable for unmet hot water load [kW] over the time interval $[t,(t+1)]$.
\end{itemize}

\noindent \textbf{States:}
\begin{itemize}
\item[\textbullet] $E_{cw,t}\in \mathbb{R}_+$: Energy level of the chilled water energy storage [kWh] at time $t$.
\item[\textbullet] $E_{hw,t}\in \mathbb{R}_+$: Energy level of the hot water energy storage [kWh] at time $t$.
\end{itemize}

\bibliography{ddip_citations}           

\end{document}